  \documentclass[12pt]{amsart}
  \usepackage{latexsym} 
  \usepackage[all]{xy}
  \usepackage{amsfonts} 
  \usepackage{amsthm} 
  \usepackage{amsmath} 
  \usepackage{amssymb}
  \usepackage{pifont}  
  \usepackage{color}
  \usepackage{hyperref}
  
\newcommand{\cC}{{C}}  
\numberwithin{equation}{section}
  \def\sw#1{{\sb{(#1)}}} 
  \def\sco#1{{\sb{[#1]}}} 
  \def\su#1{{\sp{[#1]}}}  
  \def\suc#1{{\sp{(#1)}}}

   \def\check{{\sl Check.~~}} 
   \def\endcheck{\ding{51}\medskip}
  \def\endproof{\hbox{$\sqcup$}\llap{\hbox{$\sqcap$}}\medskip} 
  \def\<{{\langle}} 
  \def\>{{\rangle}}

  \def\eps{\varepsilon}

  \def\note#1{{}} 
  \def\can{{\rm can}}
  \def\canA{{\rm can}_A}
  \def\Acan{{\rm can_l}}
  \def\tcan{\widetilde{\rm can}}
   
  \def\note#1{} 
  \def\M{{\bf M}}

  \def\lrhom#1#2#3#4{{{\rm Hom}\sb{#1, #2}(#3,#4)}} 
  \def\lhom#1#2#3{{{\rm Hom}\sb{#1-}(#2,#3)}} 
  \def\rhom#1#2#3{{{\rm Hom}\sb{-#1}(#2,#3)}}
   
  \def\lend#1#2{{{\rm End}\sb{#1-}(#2)}}

  \def\lrend#1#2#3{{{\rm End}\sb{#1,#2}(#3)}}

  \def\Lhom#1#2#3{{{\rm Hom}\sp{#1-}(#2,#3)}}

  \def\can{{\rm can}}

  \def\beq{\begin{equation}} 
  \def\eeq{\end{equation}} 
  \def\DC{{\Delta_\cC}} 
  \def \eC{{\eps_\cC}}

  \def\im{{\rm Im}}

  \def\ot{{\otimes}}

  \def\te{\tilde{e}}
  
  \def\Aro{{{}^A\!\varrho}}
  \def\roA{{\varrho^A}}
  \def\wro{{}^W\!\varrho}
 \def\ch{{\rm ch}}
  \def\tw{\widetilde{w}}
   \def\tom{\widetilde{\chi}}
      \def\tI{\widetilde{I}}

 \def\ta{\tilde{a}}

   \def\tch{\widetilde{\rm ch}}
    \def\chg{{\rm chg}}
    \def\tchg{\widetilde{\rm chg}}
\def\ttau{\tilde{\tau}}

  \newcounter{zlist} 
  \newenvironment{zlist}{\begin{list}{(\arabic{zlist})}{ 
  \usecounter{zlist}\leftmargin2.5em\labelwidth2em\labelsep0.5em 
  \topsep0.6ex%\itemsep0.3ex plus0.2ex minus0.3ex 
  \parsep0.3ex plus0.2ex minus0.1ex}}{\end{list}}

  \newcounter{blist} 
  \newenvironment{blist}{\begin{list}{(\alph{blist})}{ 
  \usecounter{blist}\leftmargin2.5em\labelwidth2em\labelsep0.5em 
  \topsep0.6ex %\itemsep0.3ex plus0.2ex minus0.3ex 
  \parsep0.3ex plus0.2ex minus0.1ex}}{\end{list}} 

  \newcounter{rlist}

\def\stac#1{\raise-.2cm\hbox{$\stackrel{\displaystyle\otimes}{\scriptscriptstyle{#1}}$}}
\def\sstac#1{\otimes_{#1}}
\def\cten#1{\raise-.2cm\hbox{$\stackrel{\displaystyle\widehat{\otimes}}
{\scriptscriptstyle{#1}}$}}
\def\scten#1{\widehat{\ot}_#1}

  \def\Label#1{\label{#1}\ifmmode\llap{[#1] }\else 
  \marginpar{\smash{\hbox{\tiny [#1]}}}\fi} 
  \def\Label{\label}

  \newtheorem{proposition}{Proposition}[section] 
  \newtheorem{lemma}[proposition]{Lemma} 
   
  \newtheorem{corollary}[proposition]{Corollary} 
  \newtheorem{theorem}[proposition]{Theorem} 

  \theoremstyle{definition} 
  \newtheorem{definition}[proposition]{Definition} 
  \newtheorem{example}[proposition]{Example} 

  \theoremstyle{remark} 
  \newtheorem{remark}[proposition]{Remark}

  \newcounter{c} 
   
  \newcommand{\etyk}[1]{\vspace{-7.4mm}$$\begin{equation}\Label{#1} 
  \addtocounter{c}{1}} 
  \renewcommand{\]}{\ifnum \value{c}=1 $$\else \end{equation}\fi} 
\begin{document} 

 \title[The relative Chern-Galois character]{Strong connections 
 and the relative Chern-Galois character for corings} 
 \author{Gabriella B\"ohm}
 \address{Research Institute for Particle and Nuclear Physics, Budapest, 
 \newline\indent H-1525
 Budapest 114, P.O.B.\ 49, Hungary}
  \email{G.Bohm@rmki.kfki.hu}
  \author{Tomasz Brzezi\'nski}    
  \address{ Department of Mathematics, University of Wales Swansea, 
  Singleton Park, \newline\indent  Swansea SA2 8PP, U.K.} 
 \email{T.Brzezinski@swansea.ac.uk}   
    \date{\today}%{March 2005} 
  \begin{abstract} 
 The Chern-Galois theory is developed for corings or coalgebras 
 over non-commu\-tative
 rings. As the first step the notion of an {\em entwined extension}  as 
 an extension of algebras
 within a bijective entwining structure over a non-commutative ring is
 introduced.  
 A strong connection for an entwined extension is defined and it is shown
 to be closely related to the Galois property and to the equivariant
 projectivity 
 of the extension. A generalisation of the Doi theorem on total integrals in
 the  
 framework of entwining structures over a non-commutative ring is obtained, and
 the bearing of strong connections on properties such as faithful flatness or
 relative injectivity is revealed.  A family 
 of  morphisms between
 the $K_0$-group of the category of finitely generated projective comodules 
 of a coring  
 and  even relative cyclic homology groups of the base algebra of 
 an entwined extension with a strong connection is constructed.
 This is termed a {\em relative Chern-Galois character}. Explicit
 examples include the computation of a Chern-Galois character of 
  depth 2 Frobenius split (or separable) extensions over 
   a separable algebra $R$. 
 Finitely generated and projective modules 
  are associated to an entwined extension with a
 strong connection, the explicit form of idempotents is derived,
 the corresponding (relative) Chern characters are computed, and
 their connection with the relative Chern-Galois character is explained.
   \end{abstract} 
  \maketitle 
  \section{Introduction}
  The Chern-Connes pairing is one of the most important and powerful tools of 
  non-commutative differential geometry \cite{Con:non}. In a recent paper 
  \cite{BrzHaj:che} a formalism has been developed connecting a certain class
  of coalgebra-Galois extensions or non-commutative principal bundles 
  \cite{BrzMaj:coa} \cite{BrzHaj:coa} with the K-theoretic aspects of
  associated 
  modules or (modules of sections of) non-commutative vector bundles. The key
  tool in this formalism is the notion of a {\em strong connection} introduced
  for Hopf-Galois extensions in \cite{Haj:str} and extended to the 
  coalgebra-Galois case in \cite{BrzMaj:geo}, while the bridge between 
  coalgebra-Galois extensions and the K-theory is provided by a
  {\em Chern-Galois character}.
  
  In a series of papers \cite{Kad:D2gal}, \cite{Kad:norm} L.\ Kadison
  has shown that 
  certain depth 2 extensions of algebras fit naturally into a framework of
  generalised Galois-type extensions in which a Hopf algebra or, more
  generally, 
  a coalgebra, is replaced by a Hopf algebroid or, more generally, by
  a coring.  
  The key here is the crucial discovery made in  \cite{KadSzl:D2bgd}
  that to any depth 2 extension one can associate two 
  bialgebroids which (co)act on  the total algebra of the extension.  
  Geometrically, these new families of Galois-type extensions can be
  understood 
  as (principal) bundles in which a standard fibre is a {\em groupoid}
  (rather 
  than a group as in classical principal bundles). From the algebraic
  point of view, 
  as in the case of coalgebra-Galois extensions, these  Hopf algebroid
  extensions are given in terms of {\em Galois corings}. The Galois theory
  for Hopf algebroids has been developed and generalised further in 
  \cite{Bohm:gal} within a
  framework of  {\em entwining structures over non-commutative rings},
 introduced earlier in \cite{Bohm:int}. Galois extensions for
  Frobenius Hopf algebroids were studied in \cite{BalSzl:fgal} using
  the language of double algebras.
 
 The Chern-Galois theory is developed in \cite{BrzHaj:che} within the
 framework 
 of bijective entwining structures {\em over a field}. The aim of the
 present paper 
 is to develop the {\em relative Chern-Galois theory} (relative in the sense
 of  Hochschild) 
 within bijective entwining structures {\em over
 non-commutative rings}, and then to illustrate that the theory
 thus developed is perfectly suited for describing examples coming from
 depth 2 extensions. The outline and main results of the paper are as follows.
 
 We begin in 
  Section~\ref{sec.prelim} with a description of preliminary
  results. The topics covered include $R$-rings and corings, entwining
  structures over non-commutative rings and their entwined modules,
and the relative cyclic homology. 
  In particular we introduce the notion of
  an entwined extension in Definition~\ref{def.entw.exten} and 
   construct a new example
  of an entwining structure over a ring in Example~\ref{ex.entw} and 
  Remark~\ref{rem.ex.entw}. 
  In Section~\ref{sec.str.con} we introduce the notion of a {\em 
  strong $T$-connection} in a bijective entwining structure over a
  non-commutative 
  ring -- the first main object of studies of the present paper. 
  Theorem~\ref{thm.Tstrconn} contains the first main result of the article: it
  relates the existence of a strong connection to the Galois property and
  equivariant projectivity of an entwined extension. The notion of 
  {\em equivariant projectivity} is introduced earlier in 
  Section~\ref{sec.str.con} and is related to the existence of connections
  through 
  a generalised Cuntz-Quillen Theorem~\ref{thm.Cuntz-Quillen}. We illustrate
  the  
  notion of a strong connection on a number of examples. In particular, we
  construct an example of a strong connection in a
   depth 2 Frobenius  split (or, dually, separable) extension over
   a separable algebra $R$.
  
 Section~\ref{sec.group}  is devoted to a special case of 
 entwined extensions in which the coaction is given by a grouplike element
 in a coring.  In particular we extend Doi's theorem on total integrals
 in Proposition~\ref{lem:j-h}. We also describe how the existence of strong 
 $T$-connections is related to the direct summand property, faithful flatness and
 relative injectivity of an entwined extension. 
 
 The main results of the paper are contained in 
 Section~\ref{sec.Chern}. First in 
 Theorem~\ref{thm.galois.char}, for any $T$-flat entwined extension 
 $B\subseteq A$
 with a strong
 $T$-connection $\ell_T$ and a finitely generated and projective 
 comodule of the structure coring $\cC$,  
 we construct a family of even cycles in the $T$-relative 
 cyclic complex of $B$. These cycles give rise to a family of maps from
 the Grothendieck group of isomorphism classes of finitely generated and
 projective $\cC$-comodules to even $T$-relative cyclic homology of $B$.
 This family that a priori depends on the choice of a strong connection, is
 termed  
 a {\em relative Chern-Galois character}. We compute explicit forms of 
 the relative Chern-Galois character for cleft Hopf algebroid extensions and 
 for  depth 2 Frobenius split extensions. We then prove that,
 with additional assumptions, the value of the relative Chern-Galois character
 is equal 
 to the value of the relative Chern character of $B$. This is done by 
 associating a finitely generated (relative) projective $B$-module to an
 entwined extension with a strong $T$-connection in Theorem~\ref{thm.proj}, and
 calculating its idempotent in Theorem~\ref{thm.idem} and the corresponding
 relative Chern character in Lemma~\ref{lemma.ch}. In this case the
 relative Chern-Galois character does not depend on the choice of a
 strong $T$-connection by the cyclic-homology arguments.
   
\section{Entwined extensions and relative cyclic homology}\label{sec.prelim} 
    
    We work over an associative commutative ring $k$ with unit. All
  algebras are assumed to be associative and with a unit. Throughout,
  $R$ denotes a $k$-algebra and $\cC$ denotes an $R$-coring with a
  coproduct $\DC$ and a counit $\eC$. $\M_R$ denotes the category of 
  right $R$-modules, $\M^\cC$ denotes the category of
  right $\cC$-comodules, etc. The notational conventions for Hom-functors are
  $\rhom R - -$ for right $R$-modules, $\lhom R - - $ for left $R$-modules,
  $\lrhom RR --$ for $R$-bimodules, and similarly but with upper indices 
  for comodules of a coring. Actions and coactions are denoted by $\varrho$
  with 
  the position of the index indicating the side on which a ring or coring acts
  or coacts,  
  for example $\varrho_M$ denotes the right action of a ring on a module $M$,
  and  
  $\Aro$ denotes a left coaction of a coring on a comodule $A$ etc.
  
\subsection{$R$-rings}

An {\em $R$-ring} $A$ is a monoid in the monoidal category of
bimodules for a $k$-algebra $R$. That is, $A$ is an $R$-$R$ bimodule
equipped with an $R$-$R$ bilinear associative multiplication map
$\mu_A:A\stac{R} A\to A$ and an $R$-$R$ bilinear unit map for the
multiplication, $\eta_A:R\to A$.

Composing $\mu_A$ with the canonical epimorphism $A\stac{k} A\to
A\stac{R} A$, one obtains a $k$-algebra structure on $A$ with unit
element $1_A := \eta_A(1_R)$. Conversely, the $R$-ring structure
is determined by the $k$-algebra structure of $A$ and the unit map
$\eta_A$ via the usual coequaliser construction.

Since any module for the $k$-algebra $A$ is in particular an
$R$-module (via the algebra homomorphism $\eta_A$), the notions of
{\em modules for the $R$-ring $A$}, and for the corresponding $k$-algebra,
coincide. 

The left regular module for the base algebra $R$
extends to a left $A$-module if and only if there exists a {\em left
augmentation map} ${}_R\varrho$, i.e. a left $R$-linear retraction of $\eta_A$
satisfying ${}_R \varrho (aa')={}_R \varrho( a\eta_A(
{}_R\varrho(a')))$, for all $a,a'\in A$. 

The {\em invariants} of a left $A$-module $M$ with respect to the left
augmentation ${}_R \varrho$ are defined as elements $m\in M$ for
which $a m= \eta_A( {}_R\varrho(a)) m$ for all $a\in A$. The
$k$-module of invariants in $M$ is isomorphic to ${\rm Hom}_{A-}(R,M)$.

An $R$-ring $A$ is called {\em separable} if the multiplication map
$\mu_A:A\stac{R} A\to A$ is a split epi in the category of $A$-$A$
bimodules. It is equivalent to the existence of an $R$-$R$ bilinear
map $\zeta:R\to A\stac{R} A$ such that 
\begin{equation}\label{eq:sep}
(\mu_A\stac{R} A)\circ(A\stac{R} \zeta) = (A\stac{R} \mu_A)\circ
(\zeta\stac{R} A), \qquad
\mu_A\circ \zeta =\eta_A.
\end{equation}
The element $\zeta(1_R)$ of $A\stac{R} A$ is called a {\em
separability idempotent}. If $A$ is a separable $R$-ring such that the
unit map $\eta_A$ is injective, then $A$ is called a
separable extension of $R$.

\subsection{$R$-corings} \label{ss:cor}
An {\em $R$-coring} $\cC$ is a comonoid in the monoidal category of
bimodules for a $k$-algebra $R$. That is, $\cC$ is an $R$-$R$ bimodule
equipped with an $R$-$R$ bilinear coassociative comultiplication
map $\Delta_\cC:\cC\to \cC\stac{R} \cC$ and an $R$-$R$ bilinear counit 
map $\eC: \cC\to R$.

A {\em right comodule for an $R$-coring $\cC$} is a pair $(M,\varrho^M)$
consisting of a right $R$-module $M$ and a right $R$-linear
coassociative and counital coaction $\varrho^M:M\to M\stac{R} \cC$ of
$\cC$ on $M$. The morphisms of right $\cC$-comodules 
are right $R$-linear maps 
which are compatible with the coactions. In accordance with
our general conventions, the category of right $\cC$-comodules
is denoted by $\M^\cC$.
Left comodules for $R$-corings and their morphims are defined
analogously, and their category is denoted by ${}^\cC\M$. 

In explicit calculations we often use Sweedler's notation for coproduct
$\DC(c) = \sum c\sw 1\stac{R}c\sw 2$, right coaction 
$\varrho^M(m) = \sum m\sco 0\stac{R} m\sco 1\in M\stac{R}\cC$, and left
coaction  
${}^N\!\varrho(n) = \sum n\sw{-1}\stac{R}n\sw{0}\in \cC\stac{R} N$. 
We make notational difference between left and right coactions (note different
types of brackets used), since later on we deal with $R$-$R$ bimodules which
are 
left and right $\cC$-comodules, but not necessarily $\cC$-$\cC$ bicomodules.

   Given a $k$-algebra $B$, ${}_B\M^\cC$ denotes
  the category of $B$-$R$ bimodules that are at the same time right
  $\cC$-comodules with a left $B$-linear coaction. For any $M\in
  {}_B\M^\cC$, the right $R$-module $B\stac{k} M$ is a right
  $\cC$-comodule with the natural coaction $\varrho^{B\sstac{k} M} =
  B\stac{k} \varrho^M$, and the left $B$-multiplication ${}_M\varrho:
  B\stac{k} M\to M$ is a right $\cC$-comodule map.

By \cite[Lemma 5.1]{Brz:str}, 
the right regular module for the base algebra $R$
extends to a right $\cC$-comodule if and only if the left regular $R$-module
extends to a left $\cC$-comodule. These properties are equivalent to the
existence of a 
grouplike element in $\cC$, that is an element $e$ such that
$\Delta_\cC(e)=e\stac{R} e$ and $\eC(e)=1_R$. In terms of a
grouplike element $e$ the right and left coactions are given by
$$\varrho^R:R\to \cC, \;\; r \mapsto er,\quad \textrm{and} \quad 
{}^R\varrho:R\to \cC, \;\; r\mapsto re .
$$
The {\em coinvariants} of a right $\cC$-comodule $M$ with respect to a
grouplike element $e\in \cC$ are defined as elements $m\in M$ such that
 $\varrho^M(m)= m\stac{R} e$. The $k$-module of coinvariants in $M$ 
is isomorphic to ${\rm Hom}^{-\cC}(R,M)$ \cite[28.4]{BrzWis:cor}.
The coinvariants of left $\cC$-comodules with respect to a grouplike element
are defined analogously. In particular, the coinvariants of $R$ as a
left, and as a right $\cC$-comodule coincide, and form a unital
$k$-subalgebra $B$ of $R$ \cite[28.5 (1)]{BrzWis:cor}.
An $R$-coring $\cC$ with a grouplike element $e$ is called a {\em Galois
coring}, provided the {\em canonical map}
$$
\can_R : R\stac{B}R \to \cC, \qquad r\stac{B}r'\mapsto rer'
$$
is bijective. In this case $\can_R$ establishes an $R$-coring isomorphism
between $\cC$ and the canonical Sweedler coring $R\stac{B}R$.

For an $R$-coring $\cC$, the left $R$-dual ${^*\cC}\colon
={\rm Hom}_{R-} (\cC,R)$ is an $R$-ring with multiplication 
$$
(ff')(c)=\sum f' ( {c_{(1)}}  f({c_{(2)}}) ),\qquad
\textrm{for all\ } f,f' \in {}^*\cC,\ c\in \cC.
$$
The unit map $R\to {}^*\cC$ is given by
$
r\ \mapsto\ \left[\ c \ \mapsto \ \eC(c r)\ \right].
$
Any right
$\cC$-comodule $M$ with coaction 
$\varrho^M(m) =\sum m\sco 0\stac{R} m\sco 1$ 
determines a right module for ${}^*\cC$ via
$$ 
m f = \sum m\sco 0f(m\sco 1),\qquad 
\textrm{for all\ } m\in M, \ f\in {}^*\cC.
$$ 
Since $\cC$-colinear maps are ${}^*\cC$-linear, there is a faithful
functor $\M^\cC\to \M_{{}^*\cC}$. It is an isomorphism if and only if
$\cC$ is a finitely generated and projective left $R$-module
\cite[19.5-6]{BrzWis:cor}.

An $R$-coring $\cC$ is called a {\em coseparable coring} if the coproduct
$\Delta_\cC:\cC\to \cC\stac{R} \cC$ is a split monic in the category of
$\cC$-$\cC$ bicomodules. By \cite[26.1]{BrzWis:cor}, this is equivalent to
the existence of a {\em 
cointegral}, i.e.\ an $R$-$R$ bilinear map $\delta:\cC\stac{R} \cC\to R$, such
that 
$$
(\cC\stac{R} \delta)\circ (\Delta_\cC\stac{R} \cC)=
(\delta\stac{R} \cC)\circ (\cC\stac{R} \Delta_\cC), \qquad
\delta\circ \Delta_\cC =\eC.
$$
 
\subsection{$R$-entwining structures and entwined modules}
\label{ss:entwstr}

A {\em right entwining structure} $(A,\cC,\psi)_R$ over a $k$-algebra
$R$ consists of an $R$-ring $A$, an $R$-coring $\cC$ and an $R$-$R$
bilinear map $\psi: \cC\stac{R} A\to A\stac{R} \cC$ satisfying the
following conditions
\begin{eqnarray*}
\psi\circ (\cC\stac{R} \mu_A)&=&(\mu_A\stac{R}
\cC)\circ(A\stac{R}\psi)\circ (\psi \stac{R} A),\\
\psi\circ(\cC\stac{R} \eta_A)&=&\eta_A\stac{R} \cC,\\
(A\stac{R} \Delta_\cC)\circ \psi &=& (\psi\stac{R} \cC) \circ
(\cC\stac{R} \psi)\circ (\Delta_\cC\stac{R} A),\\
(A\stac{R} \eC)\circ \psi&=&\eC \stac{R} A.
\end{eqnarray*}
A {\em right entwined module} for a right entwining structure
$(A,\cC,\psi)_R$ is a right $A$-module and a right
$\cC$-comodule $M$ such that
$$
\varrho^M\circ \varrho_M=(\varrho_M\stac{R}\cC)\circ (M\stac{R}\psi)\circ
(\varrho^M\stac{R} A).
$$
The morphisms of entwined modules are right $A$-linear right
$\cC$-colinear maps. The category of right entwined modules for a
right entwining structure $(A,\cC,\psi)_R$ is denoted by
$\M^\cC_A(\psi)$.

Applying
the construction in \cite[Example~4.5]{Bohm:int}, 
to a right entwining structure $(A,\cC,\psi)_R$ one can associate an
$A$-coring $(A\stac{R} \cC)_\psi$ as follows. $(A\stac{R} \cC)_\psi =
A\stac{R} \cC$ with the obvious left $A$-module structure.
 The right $A$-module structure is given by
$$
(a\stac{R} c) a'=a \psi(c\stac{R} a'),\qquad \textrm{for all\ }
a\stac{R} c\in A\stac{R}C,\; a'\in A.
$$
The coproduct is $A\stac{R} \DC$ and the counit is $A\stac{R}
\eC$. Using the same line of argument as
in  \cite[Proposition~2.2]{Brz:str}, one checks that right comodules
for the $A$-coring $(A\stac{R} \cC)_\psi$ can be identified with 
right entwined modules for the entwining structure
$(A,\cC,\psi)_R$ via the identification of $M\stac{A} (A\stac{R} \cC)$
with $M\stac{R} \cC$ for any object $M$ in $\M^{(A\stac{R} \cC)_\psi}$.

Conversely, if $A\stac{R} \cC$ is an $A$-coring with obvious left
multiplication, coproduct $A\stac{R}\DC$ and counit $A\stac{R}\eC$,
then the map
\begin{equation}\label{cor.ent}
\psi: \cC\stac{R}A\to A\stac{R}\cC, \qquad c\stac{R} a\mapsto (1_A\stac{R}c)a
\end{equation}
gives rise to a right entwining structure $(A,\cC,\psi)_R$ over $R$, provided
that
\begin{equation}\label{cond.cor.ent}
\varrho_{A\stac{R}\cC}\circ (A\stac{R}\cC\stac{k}\eta_A) = 
A\stac{R}\varrho_\cC,
\end{equation}
where $\varrho_{A\stac{R}\cC}$ is the right $A$-product in $A\stac{R}\cC$
and $\varrho_\cC$ is the right $R$-product in $\cC$. The condition
\eqref{cond.cor.ent} is necessary and sufficient to ensure
 that the map $\psi$ in \eqref{cor.ent} is right
$R$-linear.

If $A$ is an entwined module for  a right entwining structure $(A,C,\psi)_R$
with the right regular action and a $C$-coaction $\roA$, 
then $\roA(1_A)$ is a grouplike element in 
$(A\stac{R} \cC)_\psi$. For an entwined module $M$, the coinvariants with 
respect to $\roA(1_A)$ are denoted by $M^{coC}$.

{\em Left entwining structures}  over a $k$-algebra
$R$ are defined in a symmetric way. A left entwining structure
${}_R(A,\cC,\psi)$ consists of an $R$-ring $A$, an $R$-coring 
$\cC$ and an $R$-$R$ bilinear map $\psi: A \stac{R} \cC\to \cC\stac{R}
A$ satisfying compatibility conditions, analogous to the right entwining case. 
Analogously, one defines the category ${}^\cC_A\M(\psi)$ of
left entwined modules. To a left entwining structure ${}_R(A,\cC,\psi)$ 
one can associate an
$A$-coring $(\cC\stac{R} A)_\psi$,  whose left comodules  can
be identified with left entwined modules. In case $A$ is itself a left
entwined module with the left regular action, 
the coinvariants of any left entwined module $M$ with respect
to the grouplike element $\Aro(1_A)\in(\cC\stac{R} A)_\psi$ are denoted by 
${}^{coC}M$.

An example of a left entwining structure over $R$ can be constructed by
the straightforward generalisation of  \cite[Theorem~2.7]{BrzHaj:coa}.
\begin{example}\label{ex.left.canon}
Consider an $R$-coring $D$, and 
an $R$-ring $A$ that is also
a left $D$-comodule with the coaction $\Aro$. Let
$$
B = \{b\in A\; |\; \forall a\in A, \; {}^A\!\varrho(ab) =
{}^A\!\varrho(a)b \} 
$$
and suppose that the (left canonical) map
$$
\Acan: A\stac{B} A\to D\stac{R} A, \qquad a\stac{B}a'\mapsto \Aro(a)a'
$$
is bijective. 
Then ${}_R(A,D,\psi_{can})$ is a left entwining structure with the  
entwining map
\begin{equation}\label{can.psil}
\psi_{can} : A\stac{R}D\to D\stac{R}A, \qquad a\stac{R}d\mapsto
\Acan(a\,\Acan^{-1}(d\stac{R}1_A)). 
\end{equation}
The corresponding coring $(D\stac{R} A)_{\psi_{can}}$ is a Galois $A$-coring
with 
respect to the grouplike element $\Aro(1_A)$, hence $A$ is a left entwined
module over ${}_R(A,D,\psi_{can})$, and $B={}^{coD}\!\!A$.
\end{example}

A right entwining structure $(A,\cC,\psi)_R$ is said to be
bijective  if $\psi$ is a
bijection.  
In this
case ${}_R(A,\cC,\psi^{-1})$ is a left entwining structure and $\psi$
is an isomorphism of corings $(\cC\stac{R} A)_{\psi^{-1}}\to (A\stac{R}
\cC)_\psi$. 

\subsection{Entwined extensions}\label{sec.entw.ext}

In view of Sections~\ref{ss:cor} and \ref{ss:entwstr},
the following statements about a bijective right entwining structure
$(A,\cC,\psi)_R$  are equivalent. 
\begin{blist}
\item $A\in\M^\cC_A(\psi)$ with the right regular $A$-action
 and the coaction $\roA:A\to A\stac{R} \cC$;
\item $A$ is a right comodule for the $A$-coring $(A\stac{R} \cC)_\psi$  with
 the right regular $A$-action and the coaction $\roA:A\to A\stac{A}
 (A\stac{R}\cC)\simeq A\stac{R}\cC$; 
\item $\roA(1_A)$ is a grouplike element in the $A$-coring $(A\stac{R}
  \cC)_\psi$;
\item $A$ is a left comodule for the $A$-coring $(A\stac{R} \cC)_\psi$  with
 the left regular $A$-action and the coaction 
$$
A\to (A\stac{R}\cC)\stac{A} A \simeq A\stac{R}\cC
\qquad a\mapsto a\roA(1_A);
$$
\item $\psi^{-1}\big(\roA(1_A)\big)$ is a grouplike element in the $A$-coring
  $(\cC \stac{R} A)_{\psi^{-1}}$;
\item $A$ is a right comodule for the $A$-coring $(\cC \stac{R}
  A)_{\psi^{-1}}$ with the right regular $A$-action and the coaction 
$$
A\to A\stac{A}(C\stac{R} A)\simeq C\stac{R} A\qquad
a\mapsto \psi^{-1}\big(\roA(1_A) \big)a \equiv \psi^{-1}\big(\roA(a)\big);
$$
\item $A$ is a left comodule for the $A$-coring $(\cC \stac{R}
  A)_{\psi^{-1}}$ with the left regular $A$-action and the coaction 
$$
A\to (C\stac{R} A)\stac{A} A\simeq C\stac{R} A\qquad
a\mapsto \psi^{-1}\big( a\roA(1_A)\big);
$$
\item $A\in {}^\cC_A\M(\psi^{-1})$ with the left regular $A$-action and
the coaction  
\begin{equation}\label{eq:Lentw}
{}^A \varrho: A\to C\stac{R} A\qquad
a\mapsto \psi^{-1}\big( a \roA(1_A)\big).
\end{equation}
\end{blist}
Suppose that $(A,\cC,\psi)_R$ is a bijective right entwining structure such
that the equivalent conditions (a)-(h) hold.
The coinvariants of $A$ as a right entwined module for
the right entwining structure $(A,\cC,\psi)_R$, and as a left
entwined module for the left entwining structure
${}_R(A,\cC,\psi^{-1})$  are the elements of the same $k$-subalgebra
of $A$,
\begin{eqnarray*}
B \!\!\!&=& \!\!\! 
A^{co\cC} =
\{b\in A\; |\; \varrho^A(b) = b\varrho^A(1_A)\} =
\{b\in A\; |\; \forall a\in A,\; \varrho^A(ba) = b\varrho^A(a)\}\\
&=&\!\!\!  
{}^{co\cC}\! A 
=
\{ b\in A \; |\; {}^{A}\!\varrho(b) =
{}^{A}\!\varrho(1_A)b\} =
\{b\in A\; |\; \forall a\in A,\; {}^A\!\varrho(ab) = {}^A\!\varrho(a)b \}.
 \end{eqnarray*}
\begin{definition}\label{def.entw.exten}
Let $(A,\cC,\psi)_R$ be a bijective right entwining structure. An algebra
extension  
$B\subseteq A$  is called  an {\em entwined extension} provided the 
equivalent conditions (a)-(h) are satisfied and $B$ is the $k$-subalgebra 
of coinvariants.
\end{definition}

If, in a bijective right entwining structure $(A,\cC,\psi)_R$ over $R$, the
coring $\cC$ 
possesses a grouplike element $e$, then $1_A\stac{R} e$ is a grouplike element
in the $A$-coring $(A\stac{R} \cC)_\psi$. Hence the conditions (a)-(h) hold. In
particular, the right $(A,\cC,\psi)_R$, and the left
${}_R(A,\cC,\psi^{-1})$-coactions on $A$ come out as, for $a\in A$,
\begin{equation} \label{eq:roA} 
\roA (a) = \psi (e\stac{R} a)\quad \textrm{and}
\quad
\Aro (a) = \psi^{-1}(a\stac{R} e). 
\end{equation}
The subalgebra $B=A^{co\cC}$ of coinvariants in $A$
comes out as
$$
B=\{\ b\in A\ \vert\ \roA(b)=b\stac{R} e\ \}\equiv
\{\ b\in A\ \vert\ \Aro(b)=e\stac{R} b\ \}. 
$$
In this situation we say that the {\em entwined extension $B\subseteq A$
is given by a grouplike element $e\in \cC$}.
\smallskip

   The following example shows that to an 
extension of algebras $B\subseteq A$ of the type 
described in Example~\ref{ex.left.canon}, 
   one can associate not only a left entwining structure over $R$, 
   but also a right entwining 
    structure over $B$, provided that the extension is faithfully flat. 
  \begin{example}\label{ex.entw}
  Let $R$ be a $k$-algebra, $D$  and  $B\subseteq A$   
  as in Example~\ref{ex.left.canon}.
  For all $d\in D$, write $\varpi(d) = \Acan^{-1}(d\stac{R}1_A)$.
  Notice that ${\varpi}$ is left $R$-linear by the left $R$-linearity
  of $\Aro$ and it is also right $R$-linear since $A$ is an $R$-ring and
  $\Acan$ is right $A$-linear.
  Suppose that $A$ is a faithfully flat right $B$-module. Let $\cC$ be
  the associated 
  Ehresmann $B$-coring (cf.\ \cite[34.13]{BrzWis:cor}), i.e.\ $\cC$ is a
  $B$-$B$ subbimodule of $A\stac{R}A$,
  $$
  \cC = \{ \sum_i a^i\stac{R}\ta^i \in A\stac{R}A \; |\; \sum_i  
  a^i\varpi(\ta^i\sw{-1})\stac{R}\ta^i\sw 0 = 
  \sum_i 1_A\stac{B}a^i\stac{R}\ta^i\},
  $$
  with the coproduct and the counit,
  $$
  \DC(\sum_i a^i\stac{R}\ta^i)= \sum_i  
  a^i\stac{R}\varpi(\ta^i\sw{-1})\stac{R}\ta^i\sw 0,
  \qquad
  \eC (\sum_i a^i\stac{R}\ta^i) = \sum_i a^i\ta^i.
  $$
  Then $(A,\cC,\psi)_B$ is a right entwining structure over $B$, where
  $$
  \psi: \cC\stac{B}A\to A\stac{B}\cC, \qquad \sum_i a^i\stac{R}\ta^i \stac{B}
  a \mapsto \sum_i  
  a^i\varpi((\ta^ia)\sw{-1})\stac{R}(\ta^ia)\sw 0,
  $$
  is the entwining map.
  \end{example}
  \check
Recall from \cite[34.13]{BrzWis:cor} that $A\stac{R} A$ is a left
entwined module for the left entwining structure
${}_R(A,D,\psi_{can})$, constructed in Example
  \ref{ex.left.canon}. Equivalently, it is a left 
comodule for the Galois $A$-coring $(D\stac{R} A)_{\psi_{can}}$ 
with coaction  
$a\stac{R}\ta\mapsto $ \break
$\sum \psi_{can}(a\stac{R} \ta_{(-1)}) \stac{R}
\ta_{(0)}$. Note that $\cC= {}^{coD}(A\stac{R} A)$.

Since $A$ is a faithfully flat right $B$-module, the functors  
$$
{}^{coD}(-): 
{}^{D}_{A}\M({\psi_{can}})\to {}_B\M\quad \textrm{and}\quad  
A\stac{B}- : {}_B\M \to {}^{D}_{A}\M({\psi_{can}})
$$
are inverse equivalences by the (left-comodule version of) the Galois
coring structure theorem \cite[Theorem~5.6]{Brz:str}.
Hence the unit of the adjunction
  $$\theta_{A\stac{R} A}: A\stac{R}A \to A\stac{B}\cC, \qquad
a\stac{R}\ta\mapsto 
  \sum a\varpi(\ta\sw{-1})\stac{R}\ta\sw 0,
  $$
  is an isomorphism of $A$-$B$ bimodules with the inverse
  $$
  \theta_{A\stac{R} A}^{-1}: A\stac{B}\cC\to A\stac{R}A, \qquad \sum_i
  a \stac{B}a^i\stac{R}\ta^i  
  \mapsto \sum_i a a^i\stac{R}\ta^i.
  $$
  Note that the map $ \theta_{A\stac{R} A}$ is a right $B$-module map by
  the definition of $B$ as left $D$-coinvariants in $A$.
  The bimodule $A\stac{R}A$ has the canonical Sweedler $A$-coring structure.
  Using this structure and the isomorphism
  $\theta_{A\stac{R} A}$, one can induce a unique $A$-coring structure on 
  $A\stac{B}\cC$, such that  $\theta_{A\stac{R} A}$ is an isomorphism of
  $A$-corings. As a result of this one finds that the coproduct and counit in
  $A\stac{B}\cC$ come out as
  $$
  \Delta_{A\stac{B}\cC} := ( \theta_{A\stac{R} A}\ \stac{B}\
  \theta_{A\stac{R} A}) 
  \circ \Delta_{A\stac{R} A} \circ  \theta_{A\stac{R} A}^{-1} = A\stac{B}\DC
  $$
  and
  $$
  \eps_{A\stac{B}\cC} :=  \eps_{A\stac{R} A} \circ  \theta_{A\stac{R} A}^{-1}
   = A\stac{B}\eC.
   $$
   Note that the fact that $ \theta_{A\stac{R} A}$ is a right $B$-module map
   implies 
   that the compatibility condition \eqref{cond.cor.ent} between the
   original and induced $B$-actions on $A\stac{B}\cC$ is satisfied. 
   The computed form of the coproduct and counit in $A\stac{B}\cC$ thus implies
   that $A\stac{B}\cC$ is a coring associated to a right
   entwining structure 
   over $B$. In particular, the entwining map $\psi$ is induced from the 
   right $A$-multiplication in $A\stac{B}\cC$ as in \eqref{cor.ent}, i.e.\
   $$\psi (c\stac{B}a) = (1_A\stac{B}c)a 
   =  \theta_{A\stac{R} A}( \theta_{A\stac{R} A}^{-1}(1_A\stac{B}c)a), 
   $$
   and 
   has the asserted form.
  \endcheck
  
  One can investigate further the entwining structure described
  in  Example~\ref{ex.entw}. The following remark summarises basic properties
  of this structure.
  
  \begin{remark}\label{rem.ex.entw}
The notation and hypotheses of Example~\ref{ex.entw} are
  assumed.
  
  (1) If the canonical left entwining structure
  ${}_R(A,D,\psi_{can})$ in  Example~\ref{ex.left.canon} is bijective, then so
  is 
  the induced right entwining structure  $(A,\cC,\psi)_B$, provided
  $A$ is a faithfully flat left $B$-module. In this case 
  $A$ is a right $D$-comodule (cf.\ equivalent conditions (a)-(h)), and the
  inverse 
  of $\psi$ comes out explicitly as
  $$
  \psi^{-1}: A\stac{B}C\to C\stac{B}A, \qquad \sum_i a\stac{B}a^i\stac{R}\ta^i
  \mapsto \sum_i(aa^i)\sco 0\stac{R}\varpi((aa^i)\sco 1)\ta^i,
$$
where $\sum a\sco 0\stac{R}a\sco 1 := \psi_{can}^{-1}(\Aro(1_A)a)$ is the 
right $D$-coaction on $A$ induced by the inverse of the entwining 
map $\psi_{can}$.

(2)  $A$ is a right entwined module for 
$(A,\cC,\psi)_B$ with the coaction
$$
\roA : A\to A\stac{B}\cC, \qquad a \mapsto \sum \varpi(a\sw{-1})\stac{R}a\sw 0,
$$
i.e.\ $\sum \varpi(1_A\sw{-1})\stac{R}1_A\sw 0$ is a grouplike element in the
associated $A$-coring $(A\stac{B}\cC)_\psi$. The coinvariants with respect to 
this element come out as
$$
T = \{s\in A\; |\; \sum s\stac{R}1_A = 1_A\stac{R}s \}.
$$
If $A$ is a faithfully flat left or right $R$-module, then $T=R$.
By the construction,  $(A\stac{B}\cC)_\psi\simeq A\stac{R}A$ as corings,
hence in this case $(A\stac{B}\cC)_\psi$ is a Galois $A$-coring.

Therefore, if $A$ is a faithfully flat left or right $R$-module, starting from
a faithfully flat left ``coring-Galois extension" $B\subseteq A$ for an
$R$-coring $D$, i.e.\ from a Galois $A$-coring
$(D\stac{R}A)_{\psi_{can}}$, we construct the associated Ehresmann
$B$-coring $C$, and conclude that $(A\stac{B}\cC)_\psi$ is a Galois
$A$-coring, i.e.\ there is a right ``coring-Galois extension" $R\subseteq A$
for the $B$-coring $\cC$. This suggests an
interesting generalisation of and, perhaps, more symmetric
approach to {\em bi-Galois extensions} (cf.\
\cite[Definition~2.6]{FisMon:fro},  
\cite[Definition~3.4]{Sch:big}). It is shown in \cite{Sch:big} that 
to any faithfully flat Hopf-Galois {\em object} one can 
associate a bi-Galois extension. This construction is restricted
to Hopf-Galois objects as only in this case the associated Ehresmann  
coring is a Hopf algebra. To obtain more general and more symmetric
situation one needs to replace coalgebras (or Hopf algebras) over
the same ring by corings over different rings. Examples of this situation,
when $A$ is a comodule algebra for two bialgebroids $D$ and $C$, are described
in terms of $R$-$B$ torsors in \cite{Hobst:PhD}.

(3) In view of (1) and (2), if the canonical left entwining 
structure
  ${}_R(A,D,\psi_{can})$ in  Example~\ref{ex.left.canon} is bijective and
  $A$ is a faithfully flat left $B$-module, then $T\subseteq A$ is an entwined
  extension in $(A,\cC,\psi)_B$.
  \end{remark}
  
A rich source of examples of entwining structures over rings and entwined
extensions is provided by  $\times_R$-bialgebras or bialgebroids introduced
in \cite{Tak:crR}  and  \cite{Lu:hgd}, and Hopf algebroids
introduced in \cite{BohmSzl:hgdax}. 
We refer to \cite[Section~31]{BrzWis:cor} and  \cite[Section 2]{KadSzl:D2bgd}
for a review of left and right bialgebroids, respectively, and to
\cite{Bohm:hgdint} for more information on Hopf algebroids.

\begin{example} \label{ex:bgdentw}
(1) Let $(H,R,s,t,\gamma,\pi)$ be a right bialgebroid. This means that
$H$ is the total algebra, $R$ is the base algebra, $s$ is the source and $t$
is the target map, $\gamma$ is the coproduct and $\pi$ is the counit. 
Let $A$ be a right comodule
  algebra with coaction $\roA(a)=\sum a_{[0]}\stac{R} a_{[1]}$. Since the unit
  of  $H$ is a grouplike element, 
  $A$ is an entwined module for the entwining structure over
  $R$ consisting of the $R$-ring $A$, the $R$-coring
  $(H,\gamma,\pi)$ and the entwining map
\begin{equation}\label{eq:entw}
\psi:H\stac{R} A\to A\stac{R} H,\qquad h\stac{R} a\mapsto \sum a\sco 0\stac{R}
h a \sco 1.
\end{equation}

(2) Recall from  \cite{Bohm:hgdint} that a  Hopf algebroid ${\mathcal H}$ 
consists of a left bialgebroid
${\mathcal H}_L=(H,L,s_L,t_L,$ $\gamma_L,\pi_L)$, right bialgebroid 
${\mathcal H}_R=(H,R,s_R,t_R,\gamma_R,\pi_R)$ and a map $S:H\to H$,
called an {\em  antipode}, satisfying a number of compatibility conditions. 
Let $A$ be a right ${\mathcal H}$-comodule algebra in the sense of 
\cite[(arXiv version), Definition 2.8]{BB:cleft}.
This means that $A$ is both a right ${\mathcal 
  H}_R$-comodule algebra and a right ${\mathcal H}_L$-comodule algebra,
such that the ${\mathcal H}_R$-coaction 
$a \mapsto \sum a_{[0]}\stac R a_{[1]}$
is ${\mathcal H}_L$-colinear and the
${\mathcal H}_L$-coaction 
$a\mapsto \sum a_{\langle 0\rangle}\stac{L}a_{\langle 1\rangle}$
is ${\mathcal H}_R$-colinear. 

Under the additional assumption that the antipode is bijective, the entwining
map (\ref{eq:entw}) was shown in \cite[Lemma 4.1]{Bohm:gal} to be bijective
with the inverse 
$$
\psi^{-1}:A\stac{R} H\to H\stac{R} A,\qquad a\stac{R} h\mapsto 
\sum h S^{-1}(a_{\langle 1\rangle})\,\stac{R}\, a_{\langle 0\rangle}.
$$
Hence in this case $A$ is an entwined extension of its 
${\mathcal H}_R$-coinvariants (which are, in turn, the same as the ${\mathcal
H}_L$-coinvariants, see 
\cite[(arXiv version), Proposition 2.4]{BB:cleft}).

(3) 
An extension $B\subseteq A$ of $k$-algebras is called a depth 2 (or
D2, for short) extension if $A\stac{B} A$ is a direct summand in a
finite direct sum of copies of $A$ both as an $A$-$B$ bimodule and 
as a $B$-$A$ bimodule. By \cite[Lemma 3.7]{KadSzl:D2bgd}, 
 this is equivalent to 
the existence of right and left {\em  D2 quasi-bases}. A right D2 quasi-basis 
consists of  finite sets 
$\{ \gamma_j\}_{j\in J}$ in ${\rm End}_{B,B}(A)$
and  
$\{\sum_{m\in M_j} {c^j}_m\stac{B} {c^{\prime j}}_m\}_{j\in J}$ in the
commutant 
$(A\stac{B} A)^B$ of $B$ in the $A$-$A$ bimodule $A\stac{B} A$,
satisfying
$$
\sum_{j\in J,{m\in M_j}} a\gamma_j(a'){c^j}_m\stac{B} {c^{\prime j}}_m
=a\stac{B} a', 
$$
for all $a\stac{B} a'\in A\stac{B} A$. A left D2
quasi-basis is defined symmetrically. Examples of D2 extensions include 
Hopf-Galois extensions with finitely generated and projective Hopf algebras, 
centrally projective extensions and H-separable extensions 
(cf.\ \cite{KadSzl:D2bgd} for more details).

It is shown in \cite[Sections 4 and 5]{KadSzl:D2bgd} that a depth 2
extension $B\subseteq A$ determines a dual pair of 
finitely generated projective bialgebroids, with total algebras $(A\stac{B}
A)^B$ (with multiplication inherited from $A^{op}\stac{k} A$) and
$\lrend{B}{B}{A}$ (with multiplication given by composition), 
respectively. The base algebra is in both cases $R$, the commutant of $B$ in
$A$. 

$A$ is a right $(A\stac{B} A)^B$-comodule algebra, and
if $A$
is a balanced right $B$-module (e.g.\ if $A$ is a right $B$-generator), then
the coinvariants coincide with $B$. The corresponding (right) entwining map,
described in (1), is given in terms of a right D2 quasi-basis as a map
$$
(A\stac{B} A)^B\stac{R} A\to A \stac{R} (A\stac{B} A)^B,\quad
\sum_n (x_n\stac{B} x'_n)\ \stac{R}\  a\mapsto
\sum_{j\in J,{m\in M_j},n} \gamma_j(a)\
\stac{R}\ ({c^j}_m  x_n\stac{B} x'_n  {c^{\prime j}}_m ).
$$
As a consequence of the D2 property, this entwining map is bijective. The
inverse has a similar form in terms of the left D2 quasi-basis.
In particular, a balanced depth 2 extension is an entwined extension. 

Similarly,
${\mathcal  E}\colon \!\!=\lend{B}{A}$ is a left $\lrend{B}{B}{A}$-comodule
algebra. The subalgebra of coinvariants consists of the multiplications with
the elements of $A$ on the right. 
This comodule algebra determines a bijective left
entwining structure.  

(4) In the case when a D2 extension $B\subseteq A$ is also
a Frobenius extension, the bialgebroids $(A\stac{B} A)^B$ and
$\lrend{B}{B}{A}$ are  Hopf algebroids
with two-sided non-degenerate integrals, and hence with bijective
antipodes (cf.\  \cite{BohmSzl:hgdax,Bohm:alt}). In this case the inverse of
the entwining map has a simple form as in part (2).
\end{example}

\subsection{Connections and relative cyclic homology}\label{sec.cyc}
 The concepts of a 
connection and cyclic homology were introduced by A.\ Connes in
\cite{Con:non}.  
The notion of a relative cyclic homology was introduced 
by L.\ Kadison in \cite{Kad:PhD,Kad:cyc}.

Given a $T$-ring $B$ with product
  $\mu_B$, one defines a  differential graded
  algebra $\Omega B$ of {\em $T$-relative differential forms on $B$} (cf.\
  \cite[Section~2]{CunQui:alg}), by  
  $$\Omega^nB =
  \Omega^1B\stac{B}\Omega^1B\stac{B}\cdots \stac{B} \Omega^1B\qquad  
  \mbox{(n-times)},
  $$ 
  where
  $\Omega^1 B = \ker\mu_B\subseteq B\stac{T} B$. The differential is  
the $T$-$T$ bilinear map
  $d: B\to \Omega^1B$, $b\mapsto 1_B\stac{T} b - b\stac{T} 1_B,$
  and is uniquely extended as a graded differential to the whole of
$\Omega B$. Let $M$ be a left $B$-module. A $k$-linear map  
  $ \nabla_T: M\to \Omega^1B\stac{B} M \simeq  \Omega^1B M$
  is called a {\em $T$-connection}, provided it satisfies the {\em Leibniz
rule}, i.e.\ for all $b\in B$ and $m\in M$, 
  $$
  \nabla_T(bm) = d(b)\stac{B}m + b\nabla_T(m).
  $$
Notice that a $T$-connection is necessarily left $T$-linear.
There is a close relationship between $T$-connections in $M$ and the
$T$-relative projectivity of $M$
revealed by the Cuntz-Quillen Theorem \cite[Proposition~8.2]{CunQui:alg}: 
a $T$-connection
in $M$ exists if and only if there is a left $B$-linear section of the
multiplication map $B\stac{T}M\to M$. 

For any $T$-ring $B$ one defines an $(n+1)$-fold {\em circular tensor product} 
$B^{\scten{T}(n+1)}$ as the $T$-$T$ bimodule $B^{\sstac{T}(n+1)}$ 
factored by the submodule generated by
$b_0\stac{T} b_1\stac{T}\ldots \stac{T} b_n t - 
t b_0\stac{T} b_1\stac{T}\ldots \stac{T} b_n$ (see
\cite[1.2.11]{Lod:cyc}). Note that 
$B^{\scten{T}1}  = B/[B,T]$. On such a circular tensor product it is possible
to define 
the {\em cyclic operators}
$$
\tau_n: B^{\scten{T}(n+1)} \to B^{\scten{T}(n+1)}, \quad 
b_0\cten{T} b_1\cten{T}\ldots \cten{T} b_n\mapsto 
(-1)^n b_n\cten{T} b_0\cten{T}\ldots \cten{T} b_{n-1}.
$$
The relative cyclic homology of a $T$-ring $B$ is then
defined as the total homology
of the following bicomplex
$$
\xymatrix{\vdots \ar[d]^{\partial_3} & \vdots\ar[d]^{-\partial'_3} &\vdots
  \ar[d]^{\partial_3}  
& \vdots\ar[d]^{-\partial'_3} & 
\\ B^{\scten{T}3}\ar[d]^{\partial_2}&
B^{\scten{T}3}\ar[d]^{-\partial'_2}\ar[l]_{\ttau_2} 
& B^{\scten{T}3}\ar[d]^{\partial_2}\ar[l]_{N_2} 
&B^{\scten{T}3}\ar[d]^{-\partial'_2}\ar[l]_{\ttau_2} & \ldots\ar[l]_{N_2}\\
B^{\scten{T}2}\ar[d]^{\partial_1}&
B^{\scten{T}2}\ar[d]^{-\partial'_1}\ar[l]_{\ttau_1} 
& B^{\scten{T}2}\ar[d]^{\partial_1}\ar[l]_{N_1} 
&B^{\scten{T}2}\ar[d]^{-\partial'_1}\ar[l]_{\ttau_1} & \ldots\ar[l]_{N_1}\\
B/[B,T]&
B/[B,T]\ar[l]_{\ttau_0} 
& B/[B,T]\ar[l]_{N_0} 
&B/[B,T]\ar[l]_{\ttau_0} &\ldots \ar[l]_{N_0}
}
$$
where $\ttau_n = B^{\scten{T}(n+1)}-\tau_n$, 
$N_n = \sum_{i=0}^{n} (\tau_n)^i$, and
$$
\partial_n'(b_0\cten{T}b_1\cten{T}\ldots\cten{T} b_n) = 
\sum_{i=0}^{n-1}
(-1)^i
b_0\cten{T}b_1\cten{T}\ldots\cten{T} b_ib_{i+1}\cten{T}\ldots\cten{T} b_n, 
$$
$$
\partial_n = \partial'_n+ (-1)^nb_nb_0\cten{T}b_1\cten{T}\ldots\cten{T}
b_{n-1}. 
$$
The relative cyclic homology of $B$ is denoted by $HC_*(B|T)$.

The (canonical) surjections $\lambda_n: B^{\sstac{k}(n+1)}\to
B^{\scten{T}(n+1)}$ 
induce the epimorphism $\lambda_*: HC_*(B)\to HC_*(B|T)$, where $HC_*(B)$
is the usual cyclic homology of the $k$-algebra $B$. If $T$ is a separable
$k$-algebra, then $\lambda_*$ is an isomorphism \cite{Kad:cyc}.

  \section{Strong connections over non-commutative rings}\label{sec.str.con}
  Throughout the paper we are dealing with left or right modules with
  a compatible additional structure (e.g., left modules which are also
  right comodules, bimodules, etc.), and we are interested in properties which
  are respected by this additional structure. Thus we are led to the  
  following 
  \begin{definition}\label{def.rel.proj} An object $M\in {}_B\M^\cC$
  is called a {\em $T$-relative right $\cC$-equivariantly projective
  left $B$-module} or a {\em $(\cC,T)$-projective left $B$-module}
  for a $k$-algebra $T$ provided that $B$ is a right $T$-module and $M$ is a 
  $T$-$R$ bimodule such that the
  $\cC$-coaction is left $T$-linear, the left $B$-action is $T$-balanced
  (i.e.\ factors through $B\stac{T} M$), and there exists a left
  $B$-module, right $\cC$-comodule section of the $B$-multiplication map
  $B\stac{T} M\to M$.   
    
  A {\em $T$-relative left $\cC$-equivariantly projective right $B$-module}
  is defined in an analogous  way. 

A $B$-$B'$ bimodule $M$ is  called a {\em
$T$-relative right $B'$-equivariantly projective left $B$-module} or a
{\em $(B',T)$-projective left $B$-module} for a $k$-algebra $T$ provided 
that $B$ is a right $T$-module and $M$ is
a $T$-$B'$ bimodule such that the left $B$-action is $T$-balanced
and there exists a 
$B$-$B'$ bilinear section of the $B$-multiplication map $B\stac{T} M\to
M$.  

A {\em  $T$-relative left $B$-equivariantly projective right $B'$-module}
is defined in a similar way. 
\end{definition}
Note that in Definition~\ref{def.rel.proj} it is not assumed that $B$
is a $T$-ring (but, 
when $B$ is a $T$-ring, then a left $B$-action on $M$  factors through 
$B\stac{T} M$). On the other hand,
if there exists a $(C,T)$-projective and faithful left 
$B$-module 
(e.g.\ a $B$-ring which is a $(C,T)$-projective left $B$-module), 
then $B$ is a 
$T$-ring: the map  $T\to B$, $t \mapsto 1_B t$ is an algebra homomorphism 
and the right action of $T$ on $B$ coincides with that induced by the above
map, i.e.\ $b t=b (1_B t)$.

It is checked in the standard way that if $B$ is a $T$-ring, then
the $(\cC,T)$-projectivity of an object $M$ in ${}_B\M^\cC$ is equivalent to
the property that for any epimorphism $p\in {\rm Hom}^{-\cC}_{B-}(X,Y)$, 
which splits in ${}_T\M^\cC$, and any morphism 
$f\in {\rm  Hom}^{-\cC}_{B-}(M,Y)$ there exists a morphism 
$g\in {\rm Hom}^{-\cC}_{B-}(M,X)$ such that $p\circ g =f$.
Analogously, the $(B',T)$-projectivity of a $B$-$B'$ bimodule $M$ is
equivalent to the property that for any epimorphism 
$p\in {\rm Hom}_{B,B'}(X,Y)$, which splits as a $T$-$B'$ bimodule map, and
every $B$-$B'$ bimodule map $f:M\to Y$ there exists a $B$-$B'$ bimodule map
$g:M\to X$ such that $p\circ g =f$.

Any object $M\in {}_B\M^\cC$ is a $B$-${}^*\cC$ bimodule, where ${}^*\cC$ is 
the left dual ring. If $M$ is a $(\cC,T)$-projective left $B$-module,
then $M$ is  
a $({}^*\cC,T)$-projective left $B$-module. Furthermore, since 
forgetful functors ${}_B\M^\cC\to {}_B\M_R$ and 
${}_B\M_R\to {}_B\M$ 
preserve retractions, a $(\cC,T)$-projective left $B$-module is  
$(R,T)$-projective, and an $(R,T)$-projective left $B$-module 
is  $T$-relative projective. 
Clearly, a $T$-relative projective left $B$-module that is projective
as a left $T$-module is a projective left $B$-module.

Since separable functors reflect retractions, part (1) of the 
following Proposition \ref{prop:sep} (which is a straightforward
generalisation of \cite[3.2 Theorem 27]{CaeMilZhu:FrSe}) shows that
in the case of a separable $k$-algebra $R$, any $B$-$R$ bimodule, which
is a $T$-relative projective left $B$-module such that the right $R$-action is
left $T$-linear, is also $(R,T)$-projective.
Similarly, part (2) of Proposition \ref{prop:sep} (which is a
straightforward 
generalisation of \cite[Theorem 3.5]{Brz:str},\cite[26.1]{BrzWis:cor}) 
shows that in the case of a coseparable coring $\cC$, any object
$M\in {}_B\M^\cC$, which is an $(R,T)$-projective left $B$-module such that
the right $\cC$-coaction is left $T$-linear, is also $(\cC,T)$-projective. 
\begin{proposition}\label{prop:sep}
(1) The forgetful functor ${}_B\M_R\to {}_B\M$ is separable for any
$k$-algebra $B$ if and only if $R$ is a separable $k$-algebra.

(2) The forgetful functor ${}_B\M^\cC\to {}_B\M_R$ is separable for any
$k$-algebra $B$ if and only if $\cC$ is a coseparable coring.
\end{proposition}
\begin{proof}
(1) Suppose that the forgetful functor ${}_B\M^\cC\to {}_B\M_R$ is
    separable for all $k$-algebras $B$. In particular it is a
    separable functor for $B=k$, and then  $R$ is a separable $k$-
    algebra (cf.\ \cite[3.2 Theorem 27]{CaeMilZhu:FrSe}). 

Conversely, suppose that $R$ is a separable $k$-algebra and  construct
a functorial retraction $\Phi$ 
of the functorial morphism ${\rm Hom}_{B,R}(-,-)\to 
\lhom{B}{{}_B(-)}{{}_B(-)}$, $u\mapsto u$ as follows. For any $B$-$R$
bimodules $M$, $N$ and  a left $B$-module map $f:M\to N$,   
$$
\Phi(f)=\varrho_N\circ (f\stac{k} R)\circ (\varrho_M\stac{k} R)\circ
(M\stac{k} \zeta), 
$$
where  $\zeta:k\to R\stac{k} R$ is a $k$-linear map satisfying conditions
(\ref{eq:sep}). 

(2) If the forgetful functor is separable for all $k$-algebras $B$, then
it is separable for $B=k$ and hence $\cC$ is a coseparable coring by
\cite[Theorem 3.5]{Brz:str} or 
\cite[26.1]{BrzWis:cor}. 

The proof of the converse is completely analogous to the proof in 
\cite[Theorem 3.5]{Brz:str} or \cite[3.29]{BrzWis:cor}. 
Explicitly, a functorial retraction $\Phi$ of the functorial morphism ${\rm
Hom}^{-\cC}_{B-}(-,-)\to  
\lrhom{B}{R}{{}_B(-)_R}{{}_B(-)_R}$, $u\mapsto u$
is given by the map, associating to a $B$-$R$ bimodule map $f:M\to N$ 
the $B$-linear $\cC$-colinear map 
$$
\Phi(f)=(N\stac{R} \delta)\circ (\varrho^N\stac{R}\cC)\circ (f\stac{R}
C)\circ \varrho^M, 
$$
where $M$ and $N$ are objects in ${}_B\M^\cC$ and $\delta$ is a
cointegral for $\cC$.
\end{proof}

Finally, let us note that if $B$ is a separable $T$-ring, then any 
object in ${}_B\M^\cC$ is a $(\cC,T)$-projective left $B$-module. Indeed, a
left $B$-linear, right $\cC$-colinear splitting of the left $B$-multiplication
$B\stac{T} M\to M$ can be constructed in terms of a separability idempotent
$\sum_{l\in L} e_l \stac{T} f_l\in B\stac{T} B$ as $m\mapsto \sum_{l\in L}
e_l\stac{T} f_l m$.

  Following the same line of argument as in \cite{CunQui:alg}, one
  easily establishes a  
  relationship between $(\cC,T)$-projectivity and the existence of
  $T$-connections. 
For any $M\in {}_B\M^\cC$,  $ \Omega^1B\stac{B} M$ is a right
$\cC$-comodule with the natural coaction $\varrho^{ \Omega^1B\sstac{B} M}
=  \Omega^1B\stac{B} \varrho^{M}: \Omega^1B\stac{B} M\to \Omega^1B\stac{B}
M\stac{R}\cC$. 
 
 \begin{theorem}[Cuntz-Quillen] 
Let $B$ be a $T$-ring and $\cC$ an $R$-coring. An object
$M\in {}_B\M^\cC$ is a
$(\cC,T)$-projective left $B$-module if and only if there exists a
right $\cC$-colinear $T$-connection $ \nabla_T: M\to \Omega^1B\stac{B} M$. 
\label{thm.Cuntz-Quillen}
 \end{theorem}
 \begin{proof}  Note that $\Omega^1BM$ can be identified with the kernel
 of the $B$-product
 ${}_M\varrho:B\stac{T}M\to M$, by the
 commutativity of the following diagram: 
 $$\xymatrix{
  & 0\ar[d] && \\
  0\ar[r] &\Omega^1BM \ar@<.5ex>[rr]^{\iota_1}\ar[d]  &&
     B\stac{T}M\ar@{=}[d]\ar@<.5ex>[ll]^{\pi}
    \\
 0\ar[r] &\ker{}_M\varrho \ar@<.5ex>[rr]^{\iota_2}  &&
     B\stac{T}M\ar@<.5ex>[ll]^{\pi}
} $$
where $\iota_1$, $\iota_2$ are obvious inclusions both split by $\pi:
b\stac{T}m\mapsto b\stac{T} m-1\stac{T}bm.$ 

Suppose that $M$ is a $(\cC,T)$-projective left $B$-module, and let
$\sigma_T : M\to B\stac{T} M$ be a left $B$-linear, right $\cC$-colinear
section of the $B$-multiplication map ${}_M\varrho$. Then 
$$\nabla_T:M\to \Omega^1B\stac{B} M, \qquad 
m\mapsto 1_B\stac{T} m - \sigma_T(m)$$
is a right $\cC$-colinear $T$-connection. Conversely, given a $T$-connection
$\nabla_T$, define 
$$
\sigma_T: M\to B\stac{T} M, \qquad m\mapsto 1_B \stac{T} m - \nabla_T(m).
$$ 
Using the Leibniz rule for $\nabla_T$ one easily checks that $\sigma_T$ is left
$B$-linear. It obviously splits the product (as
$\nabla_T(M)\subseteq\ker{}_M\varrho$) and is right $\cC$-colinear as a
difference of right $\cC$-comodule maps. 
 \end{proof}
 
 Using arguments similar to these in the proof of
 Theorem~\ref{thm.Cuntz-Quillen} 
 one easily proves that a $B$-$B'$ bimodule $M$, where $B$ is a $T$-ring, 
 is a $(B',T)$-projective left $B$-module if and only if there exists a
 right $B'$-linear connection $ \nabla_T: M\to \Omega^1B\stac{B} M$. 

The following definition introduces the main object studied in the
present paper. 
 
\begin{definition} \label{def.str.con} 
Let $(A,\cC,\psi)_R$ be a bijective right entwining structure over 
$R$ and $B\subseteq A$ an
 entwined extension. Let $T$ be a $k$-subalgebra of $B$.
View $A\stac{T} A$ as a right $\cC$-comodule via 
$A\stac{T}\varrho^{A}$ and
as a left $\cC$-comodule via 
$ {}^{A}\!\varrho\stac{T}A$, where $\roA$ is the right and $\Aro$ is the
left $\cC$-coaction on $A$, related by \eqref{eq:Lentw}.
A left and right $\cC$-comodule map
$\ell_T:\cC\to A\stac{T} A$ is called  a {\em strong $T$-connection in
$(A,\cC,\psi)_R$} iff, for all $c\in \cC$, 
 $$
 \tcan_T(\ell_T(c)) = 1_A\stac{R} c,
 $$
where $\tcan_T: A\stac{T} A\to A\stac{R}\cC$, 
$a\stac{T}a'\mapsto a\varrho^A(a')$.
\end{definition} 
 \begin{lemma}\label{lemma.s.str}
 Let $(A,\cC,\psi)_R$ be a
bijective right entwining structure over $R$ and $B\subseteq A$ an
 entwined extension.  
Let $T'\subseteq T\subseteq B$ be $k$-subalgebras. If $A$ is a
 $(\cC,T')$-projective left or right $T$-module and there exists a
 strong $T$-connection in $(A,\cC,\psi)_R$, then there exists a
 strong $T'$-connection in $(A,\cC,\psi)_R$.
\end{lemma}
 \begin{proof}
 Suppose that $A$ is a $(\cC,T')$-projective left $T$-module, and let 
 $\xi$ be a left $T$-module right $\cC$-comodule
 splitting of the multiplication map $T\stac{T'} A\to A$. Define a
 right $\cC$-comodule map  
\begin{equation}\label{eq:ltosi}
 \ell_{T'} = (A\stac{T}\xi)\circ\ell_T\ :\ \cC\to A\stac{T'} A,
\end{equation}
where the canonical isomorphism $A\stac{T}T\to A$ is implicitly 
used (here and also in the computations below). Since
$\Aro$ is a right $T$-module map, the canonical isomorphism  
$A\stac{T}T\to A$ is a left $\cC$-comodule map. Thus $\ell_{T'}$ 
is a left $\cC$-comodule map. 
Using the right $\cC$-colinearity of $\xi$ (in the second equality)
and the fact that $\xi$ is a section of the product 
$T\stac{T'} A\to A$  (in the third), we can compute
\begin{eqnarray*}
\tcan_{T'}\circ\ell_{T'} 
&=& (\mu_A\stac{R}\cC)\circ(A\stac{T'}\roA)\circ
(A\stac{T}\xi)\circ\ell_T\\ 
&=& (\mu_A\stac{R}\cC)\circ(A\stac{T}\xi \stac{R}\cC)\circ  
(A\stac{T}\roA)\circ\ell_T\\
&=& (\mu_{A}\stac{R}\cC)\circ  
(A\stac{T}\roA)\circ\ell_T = \tcan_T\circ\ell_T,
\end{eqnarray*}
(note that $\mu_A$ denotes the product of $A$ both as a $T$- and as a 
$T'$-ring).
Thus we conclude that, for all $c\in \cC$, $\tcan_{T'}(\ell_{T'}(c)) =
1_A\stac{R} c$, i.e.\ $\ell_{T'}$
is a strong $T'$-connection in $(A,\cC,\psi)_R$. The $(\cC,T')$-projective 
 right $T$-module case is proven in a similar way. 
 \end{proof}
\begin{remark}\label{rem:Tgalcom}
In \cite[Theorem 4.4]{Brz:gal} it has been shown for a general class
of comodules of corings that the canonical map is an isomorphism (of
corings) provided that it is a split epimorphism of
comodules. Originally, \cite[Theorem 4.4]{Brz:gal} was formulated for
base algebras over fields. The proof was extended to base algebras
over commutative rings $k$ in \cite[Theorem 2.1]{BrzTur:str}. Notice,
however, that in the proof of \cite[Theorem 2.1]{BrzTur:str} the
commutativity of $k$ does not play any role. Hence, by repeating the
arguments there, one can prove the following (cf. \cite[5.9]{Wis:gal}). 

Let $D$ be an $A$-coring with a grouplike element $e$. Denote the
corresponding coinvariants of $A$ by $B$. Let $T$ be a $k$-subalgebra of
$B$. Suppose that the obvious inclusion 
$B\stac{T} A\to \Lhom D A {A\stac{T} A} \simeq 
\{\sum_i a_i\stac{T} {\tilde a}_i\in A\stac{T} A \; |\; 
\sum_i e a_i\stac{T} {\tilde a}_i = \sum_i a_ie\stac{T} {\tilde a}_i\},$
$b\stac{T} a\mapsto b \stac{T} a$,
is an isomorphism. Then D is a Galois A-coring if 
     $$
     {\widetilde {\rm can}}_T : A\stac{T} A \to D ,\quad 
     a\stac{T} a' \mapsto a \roA(a') = aea',
     $$
     is a split epimorphism of left D-comodules.
\end{remark}
In light of the bijective correspondence between $T$-connections and 
 $(\cC,T)$-projective modules described in Theorem~\ref{thm.Cuntz-Quillen}
 the following theorem, which is the main result of this section,
 justifies the use  
 of the term ``$T$-connection" in Definition~\ref{def.str.con}.

Note that, for an entwined extension $B\subseteq A$ in a bijective right
entwining structure $(A,\cC,\psi)_R$ and a subalgebra $T$ of $B$,
$A\stac{T} A$ is a right entwined module for $(A,\cC,\psi)_R$ 
with the coaction $A\stac{T}\roA$ and a  left entwined module for
${}_R(A,\cC,\psi^{-1})$ with the coaction $\Aro\stac{T}A$, where 
$\Aro$ is related to $\roA$ as in \eqref{eq:Lentw}, and with obvious
$A$-multiplications.

\begin{theorem}\label{thm.Tstrconn}
Let $(A,\cC,\psi)_R$ be a bijective right entwining structure over $R$ and let
$B\subseteq A$ be an entwined extension. Let $T$ be a $k$-subalgebra of
$B$. Consider the following statements. 
\begin{blist}
\item There exists a strong $T$-connection in $(A,\cC,\psi)_R$.
\item $(A\stac{R}\cC)_{\psi}$ is a Galois $A$-coring and $A$ is a
  $(\cC,T)$-projective left $B$-module.
\item $(\cC\stac{R}A)_{\psi^{-1}}$ is a Galois $A$-coring and $A$ is a
  $(\cC,T)$-projective right $B$-module.
\end{blist}
Then
\begin{zlist}
\item The statement (b) implies (a).
\item The statement (c) implies (a).
\item If the obvious inclusion
 $$
  B\stac{T} A \to 
  {}^{co\cC}(A\stac{T} A)
   ,\qquad 
  b\stac{T} a\mapsto b\stac{T}a, 
  $$ 
is an isomorphism, then (a) implies (b).
\item If the obvious inclusion
   $$
  A\stac{T} B \to 
  {}(A\stac{T} A)^{co\cC} , \qquad 
a\stac{T} b\mapsto  a\stac{T} b, 
  $$
is an isomorphism, then  (a) implies (c).
\end{zlist}
 \end{theorem}
 \begin{proof} 
(1) We show that if $(A\stac{R}\cC)_\psi$ is a Galois $A$-coring,
   then there exists a strong $B$-connection in the bijective right entwining
   structure $(A,\cC,\psi)_R$. Then the
   claim follows by Lemma \ref{lemma.s.str}.

If $(A\stac{R}\cC)_\psi$ is a Galois $A$-coring, then there 
exists an $A$-coring  
 inverse $\canA^{-1}$ of the canonical map $\canA: A\stac{B}A\to
 A\stac{R}\cC$,  
 $a\stac{B}a'\mapsto a\roA(a')$. Let 
 $$
 \varpi: \cC\to A\stac{B} A, \qquad c\mapsto \canA^{-1}(1_A\stac{R} c),
 $$
 be the translation map and write $\varpi(c) = \sum c\su 1\stac{B} c\su
2$. Since  
 $\canA^{-1}$ is a coring map, it is, in particular, a morphism of
right $(A\stac{R}\cC)_\psi$-comodules, hence of right
$\cC$-comodules. Therefore, also $\varpi$ 
 is a morphism of right $\cC$-comodules. Furthermore, $\canA^{-1}$ is a 
 morphism of left $A$-modules. Since $A$ is an $R$-ring, this implies
that $\varpi$ is a left $R$-module map. 
To exploit further the $A$-coring map property of $\canA^{-1}$, start
from the  
 equality
 $$
 \sum \canA^{-1}(1_A\stac{R} c\sw 1)\stac{A} \canA^{-1}(1_A\stac{R} c\sw 2) =
\sum c\su 1\stac{B} 1_A\stac{A} 1_A\stac{B} c\su 2, 
 $$
 and apply $(\psi^{-1}\stac{A} A\stac{B} A)\circ (\canA \stac{A} A\stac{B} A)$
 to arrive at the left $\cC$-colinearity of $\varpi$ (note that the left
$\cC$-coaction $\Aro$ is right $B$-linear). Thus $\varpi$ is a strong
$B$-connection, hence Lemma~\ref{lemma.s.str}
implies that there is a strong $T$-connection in $(A,\cC,\psi)_R$.

(3) Let $\ell_T: \cC\to A\stac{T} A$ be a strong $T$-connection. For all $c\in
\cC$, 
write $\ell_T(c) = \sum c\suc 1\stac{T} c\suc 2$ and define 
$$
\sigma_T: A\to A\stac{T} A, \qquad
\ a\mapsto 
\sum a\sco 0\ell_T(a\sco 1) =
\sum a\sco 0a\sco 1\suc 1\stac{T} a\sco 1\suc 2. 
$$
For all $a\in A$,
 \begin{eqnarray*}
 (\Aro\stac{T} A)\circ\sigma_T(a) 
&=& \sum  \Aro(a\sco 0a\sco 1\suc 1)\stac{T} a\sco 1\suc 2\\ 
&=& \sum \psi^{-1}(a\sco 0\stac{R} a\sco 1\suc 1\sw{-1})a\sco 1\suc
1\sw{0}\stac{T} a\sco 1\suc 2\\ 
 &=& \sum \psi^{-1}(a\sco 0\stac{R} a\sco 1\sw 1)a\sco 1\sw 2\suc 1\stac{T}
a\sco 1\sw 2\suc 2\\ 
 &=&\sum \psi^{-1}(a\sco 0\sco 0\stac{R} a\sco 0 \sco 1)a\sco 1\suc
1\stac{T} a\sco 1\suc 2\\ 
 &=& \sum \psi^{-1}(\roA(1_A))a\sco 0a\sco 1\suc 1\stac{T} a\sco 1\suc 2
=   \Aro(1_A)\sigma_T(a), 
 \end{eqnarray*}
 where the second equality follows by the fact that $A$ is a left
entwined module and the third one follows by the left colinearity of
$\ell_T$. The penultimate equality is a consequence of the 
right $A$-linearity of the right coaction $\psi^{-1}\circ \roA$ of the
$A$-coring $(\cC\stac{R} A)_{\psi^{-1}}$ on $A$.
Thus, for all $a\in A$, 
$\sigma_T(a)$ is a coinvariant of the left $(A\stac{R}\cC)_\psi$-comodule
$A\stac{T} A$,
hence, 
  $$
  \sigma_T: A\to B\stac{T} A.
  $$
  Since $\ell_T$ is a right $\cC$-colinear map, so is
$\sigma_T$. Furthermore, $\sigma_T$ is a left $B$-linear map, as, by the
definition of $B$, $\roA$ is a left $B$-linear map. The
splitting property of $\ell_T$, $ \tcan_T(\ell_T(c)) = 1_A\stac{R} c$,
 implies that, for all $c\in \cC$,  $\sum c\suc 1c\suc 2 = 1_A\eC(c)$, hence
 $$
 \mu_A(\sigma_T(a)) = \sum a\sco 0a\sco 1\suc 1a\sco 1\suc 2 = \sum
a\sco 0\eC(a\sco 1) = a. 
 $$
 Thus we conclude that $\sigma_T$ is a left $B$-linear right
$\cC$-colinear section of the multiplication map $B\stac{T} A\to A$.

 It remains to show that the right canonical map $\canA$ is an
isomorphism of $k$-modules. To this end we study first the properties
of the lifted canonical map $\tcan_T$. 
Since $A$ is a left entwined module for the left entwining structure
${}_R(A,\cC,\psi^{-1})$ with the coaction (\ref{eq:Lentw}),
$A\stac{T} A$ is a left
$(A\stac{R}\cC)_\psi$-comodule with the coaction $a\stac{T}a'\mapsto \sum
a1_A\sco 0\stac{R} 1_A\sco 1\stac{T} a'$. On the other hand, $A\stac{R} \cC$ is
a left $(A\stac{R}\cC)_\psi$-comodule via the regular coaction (the
coproduct), $a\stac{R} c\mapsto \sum a\stac{R} c\sw 1\stac{A} 1_A\stac{R} c\sw
2$. Define a left $A$-module map (hence also a left $R$-module map, as
$A$ is an $R$-ring) 
 $$
 \kappa: A\stac{R} \cC\to A\stac{T} A, \qquad a\stac{R} c\mapsto a\ell_T(c).
 $$
 Since the lifted canonical map $\tcan_T$ is left $A$-linear and, for
all $c\in \cC$, $\tcan_T(\ell_T(c)) = 1_A\stac{R} c$, one immediately finds
that  
 $\kappa$ is a section of $\tcan_T$. We claim that $\kappa$ is a left
$(A\stac{R}\cC)_\psi$-comodule map. In view of the definitions of
$\kappa$, the $(A\stac{R}\cC)_\psi$-coactions and the right action of $A$
on $(A\stac{R}\cC)_\psi$, this is equivalent to that statement that for
all $a\in A$, $c\in \cC$, 
 $$
 \sum a\psi(c\sw 1 \stac{R} c\sw 2\suc 1)\stac{T} c\sw 2\suc 2= \sum ac\suc
11_A\sco 0 \stac{R} 1_A\sco 1\stac{T} c\suc 2. 
 $$
 Using the left $\cC$-colinearity of $\ell_T$ and the form (\ref{eq:Lentw}) 
of the left $\cC$-coaction $\Aro$ we can compute 
  \begin{eqnarray*}
  \sum a\psi(c\sw 1 \stac{R} c\sw 2\suc 1)\stac{T} c\sw 2\suc 2 &=& \sum
a\psi(c\suc 1\sw{-1} \stac{R} c\suc 1\sw 0)\stac{T} c\suc 2\\ 
  &=& \sum a\psi(\psi^{-1}( c\suc 11_A\sco 0 \stac{R} 1_A\sco 1))\stac{T}
c\suc 2\\ 
 &=&\sum ac\suc 11_A\sco 0 \stac{R} 1_A\sco 1\stac{T} c\suc 2,
 \end{eqnarray*}
 as required. Thus $\kappa$ is a left $(A\stac{R}\cC)_\psi$-comodule
section of $\tcan_T$. Finally note that
$$
\Lhom{(A\stac{R}\cC)_\psi} A {A\stac{T} A}\simeq  
\Lhom{(\cC\stac{R}A)_{\psi^{-1}}} A {A\stac{T} A} 
\simeq {}^{co\cC}(A\stac{T} A) = B\stac{T}A,
$$
where the last equality follows by assumption (3).
Hence, the $T$-module version of
\cite[Theorem~4.4]{Brz:gal} adapted to this situation in
Remark \ref{rem:Tgalcom} implies that $(A\stac{R}\cC)_\psi$ is a
Galois coring. 

 Assertions (2) and (4) follow by the left-right symmetry.
 \end{proof}

 The assumptions of part (3) (resp.\ (4)) in Theorem~\ref{thm.Tstrconn}
 are automatically satisfied, provided $A$ is a flat left (resp.\ right)
 $T$-module.   
 Hence Theorem~\ref{thm.Tstrconn} leads to the following
 \begin{corollary}\label{cor.con}  Let $(A,\cC,\psi)_R$ be a
bijective right entwining structure over $R$ and let $B\subseteq A$ be an
 entwined extension. Let $T$ be a $k$-subalgebra of $B$.
 If $A$ is a flat left (resp. right) $T$-module,
then the following statements are equivalent. 
 \begin{zlist}
 \item There exists a strong $T$-connection in $(A,\cC,\psi)_R$.
 \item $(A\stac{R}\cC)_\psi$ is a Galois $A$-coring and $A$ is a
$(\cC,T)$-projective left (resp.\ right) $B$-module. 
  \end{zlist}
\end{corollary}
  
Repeating the arguments in  \cite[Theorem~4.3]{Brz:gal} one proves
also 
\begin{corollary}\label{cor.con2}  Let $(A,\cC,\psi)_R$ be a
bijective right entwining structure over $R$ and let $B\subseteq A$ 
be an entwined extension. Let $T$ be a $k$-subalgebra of $B$.
Suppose that
 \begin{blist}
 \item There exists a strong $T$-connection in $(A,\cC,\psi)_R$;
 \item $A$ is a flat left (resp.\ right) $B$-module; 
 \item $A$ is a faithfully flat left (resp.\ right) $T$-module. 
 \end{blist}
Then $A$ is a faithfully flat left (resp.\ right) $B$-module.
\end{corollary}
Note that if $A$ is a projective left (resp.\ right) $T$-module, then,
 in view of Theorem~\ref{thm.Tstrconn}~(3) (resp.\ (4)), 
 the assumption (a) in Corollary~\ref{cor.con2} implies that $A$ is a
 projective  
(hence in particular flat) left (resp.\ right) $B$-module. 

We conclude this section by constructing strong connections in some
examples. 
\begin{example}\label{ex.str.con}
Consider an entwining structure $(A,\cC,\psi)_B$ constructed in
Example~\ref{ex.entw}. Assume that $A$ is  faithfully flat as  
a left $B$-module and that the canonical left entwining map $\psi_{can}$
is bijective, hence $(A,\cC,\psi)_B$ is a bijective entwining structure by
Remark~\ref{rem.ex.entw}~(1).
The canonical inclusion map
$$
\ell_R: \cC\to A\stac{R}A, \qquad \sum_i a^i\stac{R}\ta^i\mapsto 
\sum_i a^i\stac{R}\ta^i
$$
is a strong $R$-connection.
In particular, by Corollary \ref{cor.con}, if $A$ is a flat left
or right $R$-module, then $(A\stac{B}\cC)_\psi$ is a Galois coring.

Let $T$ be a $k$-subalgebra of $R$.
If $A$ is a $(\cC,T)$-projective left $R$-module, then, 
by Lemma \ref{lemma.s.str}, a strong
T-connection in $(A,\cC,\psi)_R$ can be constructed in terms of a left
$R$-linear, right $\cC$-colinear section $\sigma_T:A\to R\stac{T} A$ of the
$R$-multiplication map,
$$
\ell_T: \cC\to A\stac{T} A,\qquad \sum_i a^i\stac{R}\ta^i \mapsto 
\sum_i a^i (\eta_A\stac{T} A)\big(\sigma_T(\ta^i)\big).
$$
In particular, if $R$ is a separable $k$-algebra and 
$\sigma_k: A\to R\stac{k}A$ is
determined by a separability idempotent $\sum_{l\in L} e_l\stac{k}
f_l$, then the strong $k$-connection comes out as
$$\ell(\sum_i a^i\stac{R}\ta^i )=
\sum_{i,l} a^i\eta_A(e_l) \stac{k} \eta_A(f_l)\, \ta^i.
$$
\end{example}
\begin{example}\label{ex:cleftstr}
Let ${\mathcal H}$ be a Hopf algebroid as in Example~\ref{ex:bgdentw}~(2).
Let $A$ be a right 
${\mathcal H}$-comodule algebra
(cf.\ 
Example~\ref{ex:bgdentw}~(2)), and let $B$ be the 
${\mathcal H}_R$-coinvariants of $A$. 
The extension $B\subseteq A$ is called a 
{\em cleft ${\mathcal H}$-extension} if
\begin{blist}
\item In addition to being an $R$-ring (with unit morphism $\eta_R$), $A$
is also an $L$-ring (with unit morphism $\eta_L$) such that
  $B$ is an $L$-subring of $A$.
\item There exists a left $L$-linear right ${\mathcal H}_R$-colinear
  map $j:H\to A$ that is
  convolution invertible in the sense that there exists a left $R$-linear
  right $L$-linear map ${\tilde j}:H\to A$ such that
$$
\mu_A\circ (j\stac{R} {\tilde j})\circ \gamma_R=\eta_L\circ\pi_L, \qquad
\mu_A\circ ({\tilde j}\stac{L} {j})\circ \gamma_L=\eta_R\circ\pi_R.
$$

\end{blist}
It turns out that the ${\mathcal H}$-cleft property of an extension
$B\subseteq A$ is sufficient and necessary for the coring 
$(A\stac{R} H)_\psi$, corresponding to the entwining structure in Example
\ref{ex:bgdentw} (1), to be a Galois $A$-coring and $A$ to be isomorphic to
$B\stac{L} H$ both  
as left $B$-modules and as right ${\mathcal H}_R$-comodules 
(normal basis property). 

Let $B\subseteq A$ be a cleft extension for a Hopf algebroid
${\mathcal H}=({\mathcal H}_L, {\mathcal H}_R,S)$.
Suppose that the antipode $S$ is bijective, hence there exists a
bijective right entwining structure $(A,{\mathcal H}_R,\psi)_R$ over
$R$ as in Example~\ref{ex:bgdentw}~(1-2). Since $(A\stac{R} H)_\psi$ is
a Galois $A$-coring, Theorem~\ref{thm.Tstrconn}~(1) implies that a
sufficient condition for the existence of a strong $T$-connection in
$(A,{\mathcal H}_R,\psi)_R$ for a $k$-subalgebra $T$ of $B$ is the $({\mathcal
  H}_R,T)$-projectivity
of the left $B$-module $A$. For an ${\mathcal H}$-cleft extension
$B\subseteq A$, there is an isomorphism of $k$-modules
$$
{\rm Hom}^{-{\mathcal H}_R}_{B-}(A,B\stac{T}A)\simeq
{\rm Hom}_{L,L}(H, B\stac{T} B),
$$
hence sections of the $B$-multiplication map $B\stac{T} A\to A$ in
${}_B\M^{{\mathcal H}_R}$ are in bijective correspondence with 
$L$-$L$ bimodule maps $f_T:H\to B\stac{T} B$, such that $\mu_B\circ
f_T=\eta_L\circ \pi_L$. In terms of $f_T$, a strong $T$-connection,
given by  formula \eqref{eq:ltosi}, is
\begin{equation}\label{eq:cleftl}
\ell_T:= (\mu_A\stac{T} \mu_A)\circ ({\tilde j}\stac{L} f_T\stac{L} j) 
\circ (\gamma_L\stac{L} H)\circ \gamma_L.
\end{equation}

In particular, an $L$-$L$ bimodule map $f_L:H\to B\stac{L} B$, satisfying 
$\mu_B\circ f_L=\eta_L\circ \pi_L$ is given by
$$
h\mapsto \eta_L\big(\pi_L(h)\big)\stac{L} 1_B\equiv 1_B\stac{L} \eta_L\big(
\pi_L(h)\big).
$$
The resulting strong $L$-connection is $\ell_L=({\tilde j}\stac{L} j)\circ
\gamma_L$. 

If $L$ is a separable $k$-algebra, then any
separability idempotent 
$\sum_{l} e_l\stac{k} f_l$ determines a
strong $k$-connection via $f(h):= \sum_{l} \eta_L(\pi_L(h)e_l)
\stac{k} \eta_L(f_l)$ for $h\in H$.
\end{example}

\begin{example}\label{ex.d2.strong}
Let $B\subseteq A$ be a depth 2 balanced extension of $k$-algebras. The
$A$-coring $A\stac{R} (A\stac{B} A)^B$ and the ${\mathcal E}$-coring 
$\lrend{B}{B}{A}\stac{R} {\mathcal E}$, corresponding to the bijective
entwining
structures in Example~\ref{ex:bgdentw}~(3), were shown to be Galois corings
in \cite{Kad:D2gal}, \cite{Kad:norm}.

Let $T$ be a subalgebra of $B$ such that $A$ is a $T$-relative projective left
$B$-module. Then it follows 
by Proposition \ref{prop:sep} and Theorem \ref{thm.Tstrconn} (1) that if $R$
is a separable $k$-algebra and $(A\stac{B} A)^B$ and ${\rm End}_{B,B}(A)$ are
coseparable $R$-corings, respectively, then there exist strong $k$-connections
in the bijective entwining structures in Example \ref{ex:bgdentw} (3). As it
turns out, $(A\stac{B} A)^B$ is a coseparable $R$-coring if the D2 extension
$B\subseteq A$ is split and ${\rm End}_{B,B}(A)$ is coseparable
if the extension $B\subseteq A$ is separable.

Motivated by forthcoming Example~\ref{ex.d2.ch}, we focus on the case
when the depth 2 balanced extension $B\subseteq A$ is also a Frobenius
extension. In this case $(A\stac{B} A)^B$ is coseparable if
and only if the D2 Frobenius extension $B\subseteq A$ is split, and
$\lrend{B}{B}{A}$ is coseparable if
and only if the D2 Frobenius extension $B\subseteq A$ is
separable. Since in the latter case $A^{op} \subseteq {\mathcal E}$ 
(where the inclusion is given by the left multiplication)
is a D2
Frobenius split extension and the bialgebroids $( {\mathcal
E}\stac{A^{op}}  {\mathcal E})^{A^{op}}$ and $\lrend{B}{B}{A}$ are
anti-isomorphic, we consider the split
case only. Fix
\begin{blist}
\item a Frobenius system $\{\omega, \sum_{k\in K} u_k\stac{B} v_k\}$ for
the extension $B\subseteq A$;
\item a right D2 quasi-basis $\{\gamma_j,\sum_{m\in M_j} {c^j}_m \stac{B}
{c^{\prime j}}_m\}_{j\in J}$ for the extension $B\subseteq A$ ;
\item a $B$-$B$ bilinear splitting $\varphi$ of the extension
$B\subseteq A$;
\item a separability idempotent $\{\sum_{l\in L} e_l\stac{k} f_l\}$
for the $k$-algebra $R$.
\end{blist}
A strong $k$-connection is explicitly constructed as the map
$$
(A\stac{B} A)^{B} \to A\stac{k} A,\qquad
a\stac{B} a' \mapsto \sum_{l\in L, k\in K, j\in J,m\in M_j}
a\omega\big(\gamma_j(a')e_lu_k\big)\stac{k} \varphi 
\big(v_k f_l {c^j}_m \big) {c^{\prime j}}_m.
$$
\end{example}

\section{The relative injectivity and  grouplike elements}\label{sec.group} 
In this section we derive conditions for an entwined extension to be a
relative injective $\cC$-comodule. We begin with the following simple
generalisation of \cite[Remark~4.2]{SchSch:gen}. 
\begin{proposition}\label{prop.rel.inj}
Let $(A,\cC,\psi)_R$ be a
bijective right entwining structure over $R$ and let $B\subseteq A$ be an
entwined extension such that $(A\stac{R}\cC)_\psi$ is a Galois
$A$-coring. 
 If the extension $B\subseteq A$ splits as a right (resp.\ left)
$B$-module, then  
 $A$ is a relative injective right (resp.\ left) $\cC$-comodule.
 \end{proposition}
 \begin{proof}
Since the right $\cC$-coaction on $A$ is left $B$-linear, the
assumption that $B$ is a direct summand of $A$ as a right $B$-module
implies that $A\simeq B\stac{B} A$ is a direct summand of $A\stac{B}
A$ as a right $\cC$-comodule. By the Galois property, $A\stac{B} A$ is
isomorphic to $A\stac{R} \cC$, in particular as a right $\cC$-comodule,
hence $A$ is a direct summand of the relative injective right
$\cC$-comodule $A\stac{R} \cC$. Therefore, it is a relative injective
$\cC$-comodule. The left-sided claim is proven analogously using 
the left-right
symmetry of bijective entwining structures.
\end{proof}

The assumptions of Proposition~\ref{prop.rel.inj} take a particularly
natural form in  
the case   of entwined extensions  given by a grouplike element
$e\in \cC$. In particular, one can derive the following generalisation of 
 Doi's theorem on total integrals \cite[(1.6) Theorem]{Doi:alg}.

\begin{proposition} \label{lem:j-h}
Let $(A,\cC,\psi)_R$ be a bijective right entwining structure over $R$ and
let $e\in \cC$ be 
a grouplike element. The right $\cC$-comodule
$A$ with the coaction given by the first of equations
(\ref{eq:roA}) is $R$-relative injective if and only if there exists a
right $\cC$-comodule map $j:\cC\to A$, such that $j(e)=1_A$. The
left $\cC$-comodule with the coaction given by the second
 of equations (\ref{eq:roA}) is $R$-relative injective if and
only if there exists a left $\cC$-comodule map ${\tilde{j}}:\cC\to A$, 
such that ${\tilde{j}}(e)=1_A$.
\end{proposition}
\begin{proof}
The $R$-relative injectivity of the right $\cC$-comodule $A$ is
equivalent to the existence of a right $\cC$-colinear retraction $h$ 
of the coaction $\roA$ in  (\ref{eq:roA}), see \cite[18.18]{BrzWis:cor}. We
show that the $k$-module of  comodule maps $j$, such that $j(e)=1_A$,
is a (non-zero) direct summand in the $k$-module of $\cC$-colinear
retractions $h$ of $\roA$.

To any right $\cC$-colinear retraction $h$ of $\roA$ associate the map
\begin{equation}
j:\cC\to A,\qquad c\mapsto h(1_A\stac{R} c). \label{eq:jtoh}
\end{equation}
The map $j$  is right $\cC$-colinear by the
colinearity of 
$h$. Since $\roA(1_A)=1_A\stac{R} e$, $j$ is also normalised, as
$$
j(e)= h(1_A\stac{R} e)= h(\roA(1_A))=1_A.
$$
Conversely, given an $e$-normalised right $\cC$-comodule map $j$, define
 the map
\begin{equation}
h=\mu_A\circ (j\stac{R} A)\circ \psi^{-1}\ :\ A\stac{R}\cC\to
A. \label{eq:htoj} 
\end{equation}
The map $h$  is right $\cC$-colinear since 
\begin{eqnarray*}
\roA \circ h
&=& \roA\circ \mu_A \circ (j\stac{R} A)\circ \psi^{-1}\\
&=& (\mu_A\stac{R} C)\circ (A\stac{R} \psi)\circ (\roA\stac{R} A)\circ 
(j\stac{R} A)\circ \psi^{-1}\\
&=& (\mu_A\stac{R} C)\circ (A\stac{R} \psi)\circ (j\stac{R} C\stac{R} A)
\circ (\Delta_C\stac{R}
A)\circ \psi^{-1}\\
&=& [\mu_A\circ (j\stac{R} A)\circ \psi^{-1}\stac{R} C]\circ
(A\stac{R} \Delta_C) 
=(h\stac{R} \cC)\circ (A\stac{R}\Delta_\cC),
\end{eqnarray*}
where the second equality follows by the fact that $A$ is a right
entwined module, the third one by the right $\cC$-colinearity of $j$, 
and the penultimate one by the definition of entwining structures. Since
$j(e)= 1_A$, the map (\ref{eq:htoj}) is a retraction of $\roA$.

Starting with a normalised $\cC$-comodule map $j:C\to A$, and associating first
the map (\ref{eq:htoj}) to it  and then the map (\ref{eq:jtoh}) to the
result, we obtain the 
comodule map
$$
c\mapsto \mu_A[(j\stac{R} A)\big(\psi^{-1} (1_A\stac{R} c)\big)]
=\mu_A(j(c)\stac{R} 1_A)
=j(c),
$$
where the first equality follows from the definition of entwining
structures. 

Similar arguments apply to left $\cC$-comodule $A$ with 
the coaction $\Aro$ in
(\ref{eq:roA}). 
\end{proof}

As an immediate application of the above results, one concludes that in 
the case of a D2 extension $B\subseteq A$ in Example~\ref{ex:bgdentw}~(3), 
$\mathcal{E} = \lend BA$ is an $R$-relative injective left comodule of 
$\lrend BBA$. Indeed, recall from Example \ref{ex.d2.strong} that the 
${\mathcal E}$-coring $\lrend{B}{B}{A}\stac{R} {\mathcal E}$, corresponding 
to the bijective entwining structure in Example~\ref{ex:bgdentw}~(3), is a 
Galois coring. What is more, the extension $A^{op}\subseteq {\mathcal E}$ is
split in the category of left $A^{op}$-modules by the map 
$\alpha\mapsto \alpha(1_A)$, hence the relative injectivity of the comodule in
question follows by Proposition \ref{prop.rel.inj}.
(The obvious inclusion $\lrend BBA\subseteq \lend BA$ is a required
left colinear map $\tilde{j}$ as in Proposition \ref{lem:j-h}.)
Furthermore, if a D2 extension $B\subseteq A$ is split (by a map 
$\varphi:A\to B$, say) in the category of right $B$-modules, then 
also $A$ is a relative injective right $(A\stac{B} A)^B$-comodule (as 
$j=\varphi\stac{B} A: (A\stac{B} A)^B\to A$ is a normalised right comodule 
map).

The derivation of sufficient conditions for assumptions of
Proposition~\ref{prop.rel.inj} is the subject of the remainder of the
present section. All these conditions turn out to be closely
related to the existence of strong connections.
\begin{lemma}\label{lem:norm}
Let $(A,\cC,\psi)_R$ be a bijective right entwining structure over $R$ and
let $B\subseteq A$ be an entwined extension given by a grouplike element
$e\in \cC$. Let $T$ be a $k$-subalgebra of $B$.
Suppose that $A$ is a flat left (resp.\ right) $T$-module and 
$B$ is a flat right (resp.\ left) $T$-module.
Then
\begin{zlist}
\item If there exists a left $B$-linear, right $\cC$-colinear section
$\sigma_T$ of the multiplication map $B\stac{T} A\to A$ 
(resp.\ a right $B$-linear, left $\cC$-colinear section
${\tilde{\sigma}}_T$ of the multiplication map $A\stac{T} B\to A$), 
then $\sigma_T(1_A)\in B\stac{T} B$ 
(resp.\ ${\tilde{\sigma}}_T(1_A)\in B\stac{T} B$).
\item If there exists a strong $T$-connection $\ell_T$ in $(A,\cC,\psi)_R$,
then $\ell_T(e)\in B\stac{T} B$.
\end{zlist}
\end{lemma}
\begin{proof}
(1) Let $\sigma_T$ be a left $B$-linear, right $\cC$-colinear section of the
  multiplication map $B\stac{T} A\to A$. 
By the right $\cC$-colinearity of $\sigma_T$ and 
$\roA (1_A) = 1_A\stac{R} e$,
$$
(B\stac{T} \roA)\big(\sigma_T(1_A)\big)
=(\sigma_T\stac{R} \cC)\big(\roA(1_A)\big)
=\sigma_T(1_A)\stac{R} e,
$$
that is, $\sigma_T(1_A)\in (B\stac{T} A)^{co\cC}$. Since $B$ is a flat 
right $T$-module,  $(B\stac{T} A)^{co\cC} = B\stac{T} B$, hence 
$\sigma_T(1_A)\in B\stac{T} B$, as required. 
The claim about $\tilde{\sigma}_T$ follows by  the left-right
symmetry of bijective entwining structures.

(2) Recall from the proof of Theorem \ref{thm.Tstrconn}~(3) that, for
any strong $T$-connection $\ell_T$, there exists a section $\sigma_T\in {\rm
Hom}^{-\cC}_{B-}(A, B\stac{T} A)$ of the multiplication map such that 
$\ell_T$ is related to the translation map $\varpi$ as
$\ell_T=(A\stac{B} \sigma_T)\circ \varpi$.
Since $\roA (1_A) = 1_A\stac{R} e$, the canonical map $\canA$ is
normalised so that
$\canA(1_A\stac{B}1_A) = 1_A\stac{R} e$. Then the 
translation map is also normalised, $\varpi(e) =$
$\canA^{-1}(1\stac{R} e) = 1_A\stac{B} 1_A$. Therefore,  
$\ell_T(e)=(A\stac{B} \sigma_T)\big(\varpi(e)\big)=\sigma_T(1_A)$, which is 
an element of $B\stac{T} B$ by (1). 
\end{proof}
\begin{proposition}\label{prop:dirsum} 
Let $(A,\cC,\psi)_R$ be a bijective right entwining structure over $R$ and
let $B\subseteq A$ be an entwined extension given by a grouplike element
$e\in \cC$. Let $T$ be a $k$-subalgebra of $B$.
Suppose that
\begin{blist}
\item There exists a strong $T$-connection in $(A,\cC,\psi)_R$;
\item $A$ is a flat left (resp.\ right) $T$-module;
\item $B$ is a flat right (resp.\ left) $T$-module;
\item $B$ is a direct summand of $A$ as a left (resp.\ right) $T$-module.
\end{blist}
Then $B$ is a direct summand of $A$ as a left (resp.\ right) $B$-module. 
\end{proposition}
\begin{proof}
Consider the assumptions without parenthesis.
By Theorem \ref{thm.Tstrconn}~(3), there exists a left $B$-linear
section $\sigma_T$ of the multiplication map 
$B\stac{T} A\to A$. In terms of $\sigma_T$ and a left $T$-module splitting 
$f$ of the canonical monomorphism $B\to A$, a left $B$-linear splitting 
can be constructed as 
$$\phi\colon =\mu_B\circ (B\stac{T} f)\circ \sigma_T.$$
Since $\phi$ is a composite of left $B$-linear maps, it is left
$B$-linear, and it satisfies
$$
\phi(1_A)=\mu_B[(B\stac{T} f)\big(\sigma_T(1_A)\big)]=
\mu_A\big(\sigma_T(1_A)\big)=1_A,
$$
where the second equality follows by Lemma \ref{lem:norm}~(1).

Similarly, if the assumptions in parenthesis hold, a right $B$-linear
splitting of the canonical monomorphism 
$B\to A$ is given by ${\tilde{\phi}}=\mu_B\circ (f\stac{T}
B)\circ {\tilde{\sigma}}_T$, where ${\tilde{\sigma}}_T$ is a right $B$-linear
section of the multiplication map $A\stac{T} B\to
A$, the existence of which follows by Theorem \ref{thm.Tstrconn}~(4),
and $f$ is a right $T$-module section of the canonical inclusion.
\end{proof}

In case $A$ is a projective left (resp.\ right) $T$-module, the
$(\cC,T)$-projectivity of $A$ as a left (resp.\ right) $B$-module implies its
projectivity. This leads to the following
\begin{corollary}\label{cor:ff}
Let $(A,\cC,\psi)_R$ be a bijective right entwining structure over $R$ and
let $B\subseteq A$ be an entwined extension given by a grouplike element
$e\in \cC$. Let $T$ be a $k$-subalgebra of $B$.
Suppose that
\begin{blist}
\item there exists a strong $T$-connection in $(A,\cC,\psi)_R$;
\item $A$ is a projective left (resp.\ right) $T$-module;
\item $B$ is a flat right (resp.\ left) $T$-module;
\item $B$ is a direct summand of $A$ as a left (resp.\ right) $T$-module.
\end{blist}
Then $A$ is faithfully flat as a left (resp.\ right) $B$-module. 
\end{corollary}
\begin{proof}
By Corollary \ref{cor.con} $A$ is a $T$-relative projective
left (resp.\ right) $B$-module, hence it is projective as left 
(resp.\ right) $B$-module by assumption (b). By Proposition~\ref{prop:dirsum},
 $B$ is a direct summand of $A$ as a left (resp.\ right)
$B$-module. Then $A$ is a generator in ${}_B\M$ (resp.\ 
$\M_B$), which implies the claim.
\end{proof}
\begin{corollary}\label{prop:relinj} 
Let $(A,\cC,\psi)_R$ be a bijective right entwining structure over $R$ and
let $B\subseteq A$ be an entwined extension given by a grouplike element
$e\in \cC$. Let $T$ be a $k$-subalgebra of $B$.
Suppose that
\begin{blist}
\item there exists a strong $T$-connection in $(A,\cC,\psi)_R$;
\item $A$ is a flat left (resp.\ right) $T$-module;
\item $B$ is a flat right (resp.\ left) $T$-module;
\item $B$ is a direct summand of $A$ as a left (resp.\ right) $T$-module.
\end{blist}
Then $A$ is $R$-relative injective as a left (resp.\ right)
$\cC$-comodule.
\end{corollary}
\begin{proof}
By Proposition~\ref{prop:dirsum}, $B$ is a direct summand of $A$ 
as a left (resp.\ right) $B$-module and by Corollary \ref{cor.con},
$(A\stac{R} \cC)_\psi$ is a Galois $A$-coring. Hence the assertion follows by
Proposition~\ref{prop.rel.inj}.
\end{proof}

 \section{The relative Chern-Galois character. Associated
   modules}\label{sec.Chern} 
 The aim of this section is, given an entwined extension $B\subseteq A$ 
 with a strong $T$-connection, to associate an Abelian group map from the 
 Grothendieck group of a coring to the even $T$-relative cyclic homology of
 $B$. A family of such maps is called a {\em relative Chern-Galois character}.
 
 The Grothendieck group $K_0(\cC)$ of an $R$-coring $\cC$ is defined as
 the Abelian group of equivalence classes of left (or right) $\cC$-comodules
 that are  
 finitely generated and projective as $R$-modules. The addition is induced from
 direct sum of such comodules. 
 
 The following lemma shows that -- just as one can associate an idempotent
 matrix of elements of $B$ to a finitely generated and projective module of an
 algebra $B$ -- finitely generated and projective comodules of a coring $\cC$
 can be characterised in terms of a `coidempotent' matrix of elements of $\cC$.

 \begin{lemma}\label{lemma.e}
 Let $W$ be a left comodule of an $R$-coring $\cC$ and suppose that $W$ is
 a finitely generated and projective $R$-module. Let $\mathbf{w} =
 \{w_i\in W, \chi_i \in {}^* W\}_{i\in I}$ be a finite dual basis of
 $W$. Define 
 $$
 e_{ij} = (\cC\stac{R}\chi_j)\left[\wro (w_i)\right] = 
 \sum w_i\sw{-1}\chi_j(w_i\sw 0), \qquad i,j\in I.
 $$
 Then, for all $i,j\in I$,
 \begin{zlist}
 \item $ \wro(w_i) = \sum_{j\in I} e_{ij}\stac{R} w_j;
 $
 \item  $
 e_{ij} = \sum_{k\in I}\chi_k(w_i)e_{kj}  =   \sum_{k\in
 I}e_{ik}\chi_j(w_k); 
 $
 \item  $
 \DC(e_{ij}) = \sum_{k\in I}e_{ik}\stac{R}e_{kj}.
 $
 \end{zlist}
 \end{lemma}
 \begin{proof}
(1) By the dual basis property,
 $$
  \wro(w_i) = \sum w_i\sw{-1}\stac{R} w_i\sw 0 = 
  \sum_{j\in I} w_i\sw{-1}\chi_j(w_i\sw 0)\stac{R} w_j =  
   \sum_{j\in I} e_{ij}\stac{R} w_j.
  $$
  
  (2) Since the coaction $\wro$ is left $R$-linear, we can compute
  $$
  \wro(w_i) = \sum_{k\in I}  \wro(\chi_k(w_i)w_k) = \sum_{k\in I}
  \chi_k(w_i)w_{k}\sw{-1}\stac{R} w_k\sw 0, 
  $$
  where the first equality follows by the dual basis property. Apply
  $\cC\stac{R} \chi_j$ to obtain the first assertion. 
The equality $e_{ij} =   \sum_{k\in I}e_{ik}\chi_j(w_k)$ follows by (1).
  
  (3) The coassociativity of $\wro$ and property (1)  imply that
  $$
  \sum_{k\in I} \DC(e_{ik}) \stac{R} w_k = \sum_{k,l}
  e_{ik}\stac{R}e_{kl}\stac{R} w_l. 
  $$
  Apply $\cC\stac{R}\cC\stac{R} \chi_j$,  use the $R$-linearity of
  $\DC$ and  (2) 
  to obtain the assertion.
\end{proof}

In fact Lemma~\ref{lemma.e} has a converse, which allows one to 
reconstruct a left $\cC$-comodule $W$ that is finitely generated and
projective as an $R$-module from a finite matrix $\mathbf{e} =
(e_{ij})_{i,j\in I}$ 
of elements in $\cC$ such that $
 \DC(e_{ij}) = \sum_{k\in I}e_{ik}\stac{R}e_{kj}.
 $ Indeed, the counit property implies that the matrix 
 $\mathbf{p} = (p_{ij}=\eC(e_{ij}))_{i,j\in I}$ is an idempotent matrix, hence 
 $W = R^{(I)}\mathbf{p}$ is a finitely generated and projective left
 $R$-module. 
 The left $\cC$-coaction $\wro: W\to \cC\stac{R} W$ is then defined by
 $$
 (\sum_{i\in I}r_ip_{ij})_{j\in I} \mapsto 
 \sum_{i,k\in I}r_ie_{ik}\stac{R}(p_{kj})_{j\in I}\simeq (\sum_{i\in
 I}r_ie_{ij})_{j\in I}. 
 $$
 Thus, rather than specifying a left $\cC$-comodule $W$, we can equally well 
 specify a matrix of elements in $\cC$ that satisfy conditions 
 of Lemma~\ref{lemma.e}~(3). Note further that such a matrix also determines
 a right $\cC$-comodule that is finitely generated and projective as a right
 $R$-module. This reflects the duality between finitely generated and
 projective left and right comodules.

 \begin{definition}\label{def.Tflat}
 Let $(A,\cC,\psi)_R$ be a bijective right
 entwining structure over $R$ and let $B\subseteq A$ be an entwined extension.
 Let $T$ be a $k$-subalgebra of $B$. Define the map
 $$
 \upsilon_T: A/[A,T]\to (A\stac{R}\cC)_\psi / [(A\stac{R}\cC)_\psi,T],
 $$
 as the projection of the $T$-$T$-bimodule map
 $$A \to (A\stac{R}\cC)_\psi , \qquad 
 a \mapsto \roA(a)-a1_{[0]}\stac{R} 1_{[1]}.
 $$
 Note that, for all $b\in B$, $\upsilon_T([b]) = 0$.
 
 We say that the entwined extension $B\subseteq A$ is {\em $T$-flat} if
 $B$ and $A$ are flat as left and right $T$-modules and the map
\begin{equation}\label{ups}
 B/[B,T] \to \ker\upsilon_T, \qquad [b]\mapsto [b],
\end{equation}
 is an isomorphism.
 \end{definition}
 
 If $T=k$, then the map
 \eqref{ups}  
 is always an isomorphism. Hence, in this case, $B\subseteq A$ is $T$-flat,
 provided $A$ and $B$ are flat $T$-modules. More generally, the map
 \eqref{ups}  
 is an isomorphism if the functor
 $$
 {}_T\M_T\to \M_k, \qquad M\mapsto M/[M,T],
 $$
 preserves monics (e.g.\ if $T$ is a separable $k$-algebra).
 
 \begin{lemma}\label{lemma.B}
Let $(A,\cC,\psi)_R$ be a bijective right entwining structure over $R$ and 
let $B\subseteq A$ be an entwined extension.
Let $T$ be a $k$-subalgebra of $B$  such that  $A$ is a  flat left and
right $T$-module. For any strong
$T$-connection  
 $\ell_T: \cC\to A\stac{T} A$ in  $(A,\cC,\psi)_R$, write $\ell_T(c) = \sum
c\suc 1\stac{T} c\suc 2$, for all $c\in \cC$. Then 
 $$
\sum \ell_T(c\sw 1)\ell_T(c\sw 2) =  \sum c\sw 1\suc 1\stac{T} c\sw 1\suc
2c\sw 2\suc 1\stac{T} c\sw 2\suc 2\in A\stac{T} B\stac{T} A. 
 $$
\end{lemma}
\begin{proof}
Since $T$ is a $k$-subalgebra of $B$ and $B$ is a $k$-subalgebra of $A$, the
coaction $\roA:A\to (A\stac{R} \cC)_\psi$ being left $B$-linear is left
$T$-linear and being right $A$-linear is right $T$-linear.
Compute, for all $c\in \cC$,
\begin{eqnarray*}
(A\stac{T} \roA\stac{T} A)(\sum \ell_T(c\sw 1)\ell_T(c\sw 2) ) 
&=&  \sum c\sw 1\suc 1\stac{T} 
c\sw 1\suc 2\sco 0\psi(c\sw 1\suc 2\sco 1\stac{R}
c\sw 2\suc 1)\stac{T} c\sw 2\suc 2\\ 
&=&  \sum c\sw 1\suc 1\stac{T} c\sw 1\suc 2\psi(c\sw 2\stac{R} c\sw
3\suc 1)\stac{T} c\sw 3\suc 2\\ 
&=&  \sum c\sw 1\suc 1\stac{T} c\sw 1\suc 2\psi(c\sw 2\suc
1\sw{-1}\stac{R} c\sw 2\suc 1\sw 0)\stac{T} c\sw 2\suc 2\\ 
&=& \sum c\sw 1\suc 1\stac{T} c\sw 1\suc 2c\sw 2\suc 11_A\sco
0\stac{R} 1_A\sco 1\stac{T} c\sw 2\suc 2, 
\end{eqnarray*}
where the first equality follows by the fact that $A$ is a right
entwined module, the second and the third equalities follow by the
$\cC$-colinearity of $\ell_T$, and the final equality is a consequence
of the definition \eqref{eq:Lentw} of the left $\cC$-coaction on $A$. Since
$A$ is flat as a left and as a right $T$-module, the assertion follows. 
\end{proof}

 \begin{theorem}\label{thm.galois.char}
Let $(A,\cC,\psi)_R$ be a bijective right entwining structure over $R$ and 
let $B\subseteq A$ be a $T$-flat entwined extension. Suppose 
that there exists a strong
$T$-connection  
 $\ell_T: \cC\to A\stac{T} A$ in  $(A,\cC,\psi)_R$, and write 
 $\ell_T(c) = \sum c\suc 1\stac{T}c\suc 2$.  For any finite matrix 
 $\mathbf{e} = (e_{ij})_{i,j\in I}$ of elements of $\cC$ such that
 $ \DC(e_{ij}) = \sum_{k\in I}e_{ik}\stac{R}e_{kj}$ and for any non-negative
 integer $n$, define the following element of the circular tensor product
 \begin{eqnarray*}
\tchg^T_{n}(\mathbf{e}) &=& \sum_{i_1,\ldots, i_{n+1}\in I} \!\!\!\!
e_{i_1i_2}\suc 2\ell_T(e_{i_2i_3})\ell_T(e_{i_3i_4})\cdots 
 \ell_T(e_{i_{n}i_{n+1}})\ell_T(e_{i_{n+1}i_1})e_{i_1i_2}\suc 1\\
& \!\!\!\! \!\!\!\!=& \!\! \!\!\!\!  \!\!\!\! \!\!\!\!\sum_{i_1,\ldots,
   i_{n+1}\in I}e_{i_1i_2}\suc 
2e_{i_2i_3}\suc 1 
\cten{T} e_{i_2i_3}\suc 2e_{i_3i_4}\suc 1\cten{T}\cdots 
\cten{T} e_{i_{n}i_{n+1}}\suc 2e_{i_{n+1}i_1}\suc 1
\cten{T} 
e_{i_{n+1}i_1}\suc 2e_{i_1i_2}\suc 1.
 \end{eqnarray*}
Then there is  a family of maps of Abelian groups
$$
\chg^T_{2n} : K_0(\cC) \to HC_{2n}(B|T), \qquad [W]\mapsto 
[\bigoplus_{l=0}^{2n}(-1)^{\lfloor l/2\rfloor}
\frac{l!}{\lfloor l/2\rfloor!}\tchg_{l}^T(\mathbf{e})],
$$
where $\mathbf{e}$ corresponds to a left $\cC$-comodule $W$ as in
Lemma~\ref{lemma.e} and $\lfloor x\rfloor$ is the integer value of $x$. The
 family  $\chg^T_{2n}$ is termed a 
 {\em relative Chern-Galois character of
   the $T$-flat $(A,\cC,\psi)_R$-entwined extension $B\subseteq A$ with 
   a strong $T$-connection $\ell_T$}. 
\end{theorem}
\begin{proof}
Lemma~\ref{lemma.B} together with the assumption that $B\subseteq A$
is a $T$-flat extension implies that, for all $n\in\mathbb{N}$, 
$\tchg^T_n(\mathbf{e}) \in B^{\scten{T}(n+1)}$. Repeating the same steps as
in the proof of Lemma~\ref{lemma.B} one finds that 
$\upsilon^T([\sum_{i\in I}e_{ii}\suc 2e_{ii}\suc 1]) =0$. Since 
$B\subseteq A$
is a $T$-flat extension, this implies that 
$\tchg^T_0(\mathbf{e}) \in B/[B,T]$. 

As $\ell_T$ is a strong $T$-connection, $\sum c\suc 1c\suc 2 =
\eta_A(\eC(c))$. Note 
further that the $\tchg^T_n(\mathbf{e})$ are
 invariant under the cyclic operator (modulo sign). These two facts allow one 
to derive  the following relations (cf.\ 
Section~\ref{sec.cyc} for the definitions of the maps involved):
$$
N_n(\tchg^T_{n}(\mathbf{e})) =  (n+1)\tchg^T_{n}(\mathbf{e}), \qquad
\partial_n(\tchg^T_{n}(\mathbf{e})) =  
 \tchg^T_{n-1}(\mathbf{e})
 $$
if $n$ is even, and  
$$
\partial'_n(\tchg^T_{n}(\mathbf{e})) = 
 \tchg^T_{n-1}(\mathbf{e}), \qquad
 \ttau_n(\tchg^T_{n}(\mathbf{e}) )= 
2\tchg^T_{n}(\mathbf{e})
$$
if $n$ is odd (and 0 otherwise). These relations imply in turn that, for any
non-negative 
integer $n$,  
$\bigoplus_{l=0}^{2n}(-1)^{\lfloor l/2\rfloor}
\frac{l!}{\lfloor l/2\rfloor!}\tchg_{l}^T(\mathbf{e})$
is a cycle in the $T$-relative cyclic complex of $B$.

Next we prove that the $\tchg_n^T(\mathbf{e})$ do not 
depend on the choice of a dual basis. In addition to a basis
$\mathbf{w} = \{w_i\in W, \chi_i \in {}^* W\}_{i\in I}$ that 
leads to the matrix $\mathbf{e}$, take a
different dual basis $\mathbf{\tw} = \{\tw_k\in W, \tom_k \in
{}^* W\}_{k\in \tI}$, and let $\mathbf{\te} =(\te_{kl})_{k,l\in \tI}$ be as in
Lemma~\ref{lemma.e}, corresponding to basis $\mathbf{\tw}$. 
Then, for all $i\in
I$, 
$$
\sum_{i'\in I} e_{ii'}\stac{R} w_{i'} = 
\wro (w_i) = \sum_{k\in \tI} \wro(\tom_k(w_i)\tw_k) = \sum_{k,l\in
\tI} \tom_k(w_i) 
\te_{kl} \stac{R} \tw_l.
$$ 
Applying $\cC\stac{R} \chi_j$ and using Lemma~\ref{lemma.e}~(2), we
obtain, for all $i,j\in I$, 
$$
e_{ij} = \sum_{k,l\in \tI} \tom_k(w_i)
\te_{kl} \chi_j(\tw_l).
$$
With this relation between the $e_{ij}$ and the $\te_{kl}$ at hand,
and using the fact that a strong $T$-connection is both left and right
$R$-linear we can compute, for all $i,j\in I$, 
\begin{eqnarray*}
\sum_{p\in I}\ell_T(e_{ip})\ell_T(e_{pj}) &=&\!\!\!\!\sum_{p\in I,
 k,l,m,q\in \tI}\!\!\!\! 
 \tom_k(w_i)\te_{kl}
\suc 1\stac{T} \te_{kl}\suc 2\chi_p(\tw_l)\tom_m(w_p)\te_{mq}\suc
 1\stac{T} \te_{mq}\suc 2 
\chi_j(\tw_q)\\
&=& \!\!\!\sum_{k,l,m,q\in \tI}\!\!\!
 \tom_k(w_i)\te_{kl}
\suc 1\stac{T} \te_{kl}\suc 2\tom_m(\tw_l)\te_{mq}\suc 1\stac{T} \te_{mq}\suc 2
\chi_j(\tw_q) ,
\end{eqnarray*}
where the second equality follows by the dual basis property. Using
the left $R$-linearity of $\ell_T$ and Lemma~\ref{lemma.e}~(2) for the
$\te_{kl}$, this can be rewritten further as 
\begin{eqnarray*}
\sum_{p\in I}\ell_T(e_{ip})\ell_T(e_{pj}) &=&
\!\!\sum_{k,l,m,q\in \tI}\!\!
 \tom_k(w_i)\ell_T(\te_{kl})\ell_T(\tom_m(\tw_l)\te_{mq})
\chi_j(\tw_q) \\
&= & \!\sum_{k,l,q\in \tI}\!\tom_k(w_i)\ell_T(\te_{kl})\ell_T(\te_{lq})
\chi_j(\tw_q).
\end{eqnarray*}
Using this equality sufficiently many times we obtain, for all $n\in
\mathbb{N}$, 
\begin{eqnarray*}
\tchg_{n}^T(\mathbf{e}) &=&  \!\!\!\!\!\!\!\!\!\!\!\sum_{i\in I, i_1,\ldots,
i_{n+1},j \in \tI} \!\!\!\!\!\!\!\! 
\te_{i_1i_2}\suc 2\ell_T(\te_{i_2i_3})\ell_T(\te_{i_3i_4})\cdots 
 \ell_T(\te_{i_{n}i_{n+1}})\ell_T(\te_{i_{n+1}j})
\chi_i(\tw_{j})\tom_{i_1}(w_i) 
\te_{i_1i_2}\suc 1\\
&=&  \!\!\!\!\!\!\!\!\sum_{i_1,\ldots, i_{n+1},j \in \tI}
\!\!\!\!\!\!\!\! 
\te_{i_1i_2}\suc 2\ell_T(\te_{i_2i_3})\ell_T(\te_{i_3i_4})\cdots 
 \ell_T(\te_{i_{n}i_{n+1}})\ell_T(\te_{i_{n+1}j}\tom_{i_1}(\tw_{j}))
\te_{i_1i_2}\suc 1\\
&=& \tchg_n^T(\mathbf{\tilde{e}}),
\end{eqnarray*}
where we used the dual basis property, the right $R$-linearity of $\ell_T$ 
and Lemma~\ref{lemma.e}~(2). Similarly, for the zeroth component,
$$
\tchg_{0}^T(\mathbf{e}) = [\sum_{i\in I, j, k \in \tI} 
\te_{jk}\suc 2\chi_i(\tw_{k})\tom_{j}(w_i) 
\te_{jk}\suc 1] = 
[\sum_{j, k \in \tI} 
\te_{jk}\suc 2\tom_{j}(\tw_{k}) 
\te_{jk}\suc 1] =  \tchg_0^T(\mathbf{\tilde{e}}),
$$
where the last equality follows by the $R$-linearity of $\ell_T$ and by 
Lemma~\ref{lemma.e}~(2). Thus we conclude that the
$\tchg_n^T(\mathbf{e})$ do not
depend on the  
choice of a dual basis of $W$. 

Suppose that there is a left $\cC$-comodule isomorphism 
$\varphi: W\to \widehat{W}$. As the $\tchg_n^T(\mathbf{e})$ do
 not depend on the choice
of a dual basis, we can choose the dual basis 
$\{\widehat{w}_i = \varphi(w_i) \in \widehat{W}, 
\widehat{\chi}_i =
\chi_i\circ\varphi^{-1} \in {}^*\widehat{W} \}_{i\in
  I}$  
in $\widehat{W}$. Then
\begin{eqnarray*}
\widehat{e}_{ij} &=& 
 \sum \widehat{w}_i\sw{-1}\widehat{\chi}_j(\widehat{w}_i\sw 0) =
 \sum \varphi(w_i)\sw{-1}\chi_j(\varphi^{-1}(\varphi(w_i)\sw 0))\\
 & =&
 \sum w_i\sw{-1}\chi_j(w_i\sw 0) = e_{ij}.
\end{eqnarray*}
 The second equality follows by the left colinearity of $\varphi$. 
 Thus the components $\tchg^T_n(\mathbf{\widehat{e}})$
 of the relative Chern-Galois cycle 
 for the matrix
 $\mathbf{\widehat{e}}$ corresponding to 
 $\widehat{W}$ coincide with $\tchg^T_n(\mathbf{e})$. All
 this proves that
 the maps $\chg_{2n}^T$ are well-defined. To prove that these are group
 morphisms, 
 note that if $W = W^1\oplus W^2$,
 then the  matrix $\mathbf{e} = (e_{ij})_{i,j\in I}$ is a direct 
 sum $\mathbf{e} = \mathbf{e}^1\oplus \mathbf{e}^2$ of 
 corresponding matrices for $W^1$ and $W^2$. The 
  explicit
 form of the relative Chern-Galois cycle  then immediately implies that 
 $\tchg_n^T(\mathbf{e}) = \tchg_n^T(\mathbf{e}^1)+\tchg_n^T(\mathbf{e}^2)$,
 hence 
 $\chg_{2n}^T([ W^1\oplus W^2]) = \chg_{2n}^T([ W^1])+ \chg_{2n}^T([W^2])$, as 
 required.
\end{proof}

In case $T=k$, the components of the relative  Chern-Galois 
characters are
denoted by $\tchg_n$ and $\chg_{2n}$. We now describe  Chern-Galois
characters for examples of entwined extensions constructed in 
Section~\ref{sec.str.con}.

 \begin{example}\label{ex.str.con.ch}
Consider a strong $T$-connection $\ell_T$ in a
bijective entwining structure $(A,\cC,\psi)_B$ described
 in Example~\ref{ex.str.con}. Assume that the extension 
$A^{co\cC} \subseteq A$ is
 $T$-flat. 
  Let 
$\mathbf{e} =  (e_{ij})_{i,j\in I}$
be a finite matrix of elements of $\cC$ such that the condition (3) in
Lemma~\ref{lemma.e} is satisfied. 
Write $e_{ij} := \sum_{m\in M_{ij}} {a_{ij}}^{(m)}
\stac{R} {{\tilde a}_{ij}}~^{(m)}$ and $\sigma_T(a)=\sum a^{\{1\}}\stac{T}
a^{\{2\}}$ 
for the components of a left $R$-linear, right $C$-colinear section
of the $R$-multiplication map $R\stac{T}A\to A$ that defines $\ell_T$. 
Then
\begin{eqnarray*}
\tchg_{n}^T(\mathbf{e}) =\sum \!\!\!\!\!\!\!&&
{{\ta_{i_1i_2}}}~^{(m_1)\{2\}}{a_{i_2i_3}}^{(m_2)} 
\eta_A( \ta_{i_2 i_3}~^{(m_2)\{1\}})\cten{T} \\
&&{{\ta_{i_2i_3}}}~^{(m_2)\{2\}}{a_{i_3i_4}}^{(m_3)} 
\eta_A( \ta_{i_3 i_4}~^{(m_3)\{1\}})
\cten{T}\cdots \cten{T} \\
&&
{{\ta_{i_{2n+1}i_{1}}}}~^{(m_{n+1})\{2\}}{a_{i_{1}i_2}}^{(m_1)}
\eta_A( \ta_{i_1 i_2}~^{(m_1)\{1\}}), 
\end{eqnarray*}
where summation is over the indices: $i_1,\dots, i_{n+1}\in I$, 
$m_1\in M_{i_1 i_2},\dots, m_{n+1}\in M_{i_{n+1} i_1}$ and over the
components of the map $\sigma_T$.

In particular, if $R$ is a separable $k$-algebra and the strong
$k$-connection $\ell$ is induced by a separability idempotent $\sum_{l\in
L} e_l\stac{k} f_l$, then
\begin{eqnarray*}
\tchg_{n}(\mathbf{e}) =\sum \!\!\!\!\!\!\!&&
\eta_A(f_{l_1}){{\ta_{i_1i_2}}}~^{(m_1)}{a_{i_2i_3}}^{(m_2)} \eta_A(e_{l_2})
\stac{k}\\ 
&&\eta_A(f_{l_2}){{\ta_{i_2i_3}}}~^{(m_2)}{a_{i_3i_4}}^{(m_3)} 
\eta_A(e_{l_3})
\stac{k}\cdots \stac{k} \\
&&
\eta_A(f_{l_{n+1}}){{\ta_{i_{n+1}i_{1}}}}~^{(m_{n+1})}
{a_{i_{1}i_2}}^{(m_1)} \eta_A(e_{l_1}), 
\end{eqnarray*}
where summation is over the indices: $i_1,\dots, i_{n+1}\in I$, 
$m_1\in M_{i_1 i_2},\dots ,m_{n+1}\in M_{i_{2n+1} i_1}$ and $l_1,\dots,
l_{n+1}\in L$. 
\end{example}
\begin{example}
Let ${\mathcal H}$ be a Hopf algebroid with a bijective antipode and let 
$B\subseteq
A$ be a cleft ${\mathcal H}$-extension as in Example
\ref{ex:cleftstr}. Suppose
that there exists a strong $T$-connection (\ref{eq:cleftl}) in the bijective
entwining structure in Example~\ref{ex:bgdentw} (1-2) and assume that the
extension $B\subseteq A$ is $T$-flat. 

Let $\mathbf{e} =  (e_{ij})_{i,j\in I}$
be a finite matrix of elements of $H$ such that the condition (3) in
Lemma~\ref{lemma.e} is satisfied 
for the coproduct $\gamma_R$ of the right bialgebroid in ${\mathcal H}$, 
and  set $c\colon = \sum_{i\in I} e_{ii}$. For a map $f_T$ in Example
\ref{ex:cleftstr}, write $f_T(h)= \sum h^{(1)}\stac{T} h^{(2)}$. Then 
\begin{eqnarray*}
\tchg_{n}^T(\mathbf{e})= \!\!\!&\sum&\!\!\!
{c_{(n+3-t)}}^{(2)} {c_{(n+4-t)}}^{(1)}\cten{T}\dots \cten{T}
{c_{(n+1)}}^{(2)} {c_{(n+2)}}^{(1)}\cten{T}\\
&&\!\!\! {c_{(n+2)}}^{(2)} j(c_{(n+3)}){\tilde j}(c_{(1)}) {c_{(2)}}^{(1)}
\cten{T}
{c_{(2)}}^{(2)}{c_{(3)}}^{(1)}\cten{T}\dots \cten{T}
{c_{(n+2-t)}}^{(2)} {c_{(n+3-t)}}^{(1)},
\end{eqnarray*}
for any $t=1,\dots, n+1$, where 
$\sum c_{(1)}\stac{L}\dots \stac{L}c_{(n+3)}$ stands for the 
action of the coproduct $\gamma_L$ iterated $n+2$ times on $c$.

In the particular example of $\ell_L=({\tilde j}\stac{L} j)\circ \gamma_L$ we
obtain
$$
\tchg_{n}^L(\mathbf{e})=\sum
1_A\cten{L}\dots \cten{L} 1_A\cten{L} j(c_{(2)}){\tilde j}(c_{(1)}) \cten{L}
1_A\cten{L}\dots \cten{L} 1_A.
$$
In contrast to cleft Hopf algebra extensions, cleft Hopf algebroid extensions
provide non-trivial relative Chern-Galois characters.
\end{example}

\begin{example}\label{ex.d2.ch}
Let $B\subseteq A$ be a balanced depth 2 split Frobenius
extension of $k$-algebras such that the commutant $R$ of $B$ in $A$ is a
separable $k$-algebra. Suppose that $A$ is a flat $k$-module.
Let $W$ be a left comodule for the right bialgebroid $(A\stac{B} A)^B$
that is finitely generated and projective as a left
$R$-module and let $\mathbf{e} = (e_{ij})_{i,j\in I}$ be the corresponding
matrix of elements of $(A\stac{B} A)^B$ (cf.\ Lemma~\ref{lemma.e}). Take a
finite dual basis  
$\mathbf{w} = \{w_i\in W, \chi_i \in {}^* W\}_{i\in I}$ and introduce the
element  
$$
\sum_q {x}_q\stac{B} {x^{\prime}}_q\colon=
\sum_{i\in I} ((A\stac{B} A)^B\stac{R}
\chi^i)\,{^W\!\!\varrho}(w_i)\equiv \sum_{i\in I} e_{ii}
$$ 
of $(A\stac{B} A)^B $. Fix the data as in (a)-(d) of
Example~\ref{ex.d2.strong}. 
By Lemma \ref{lemma.e} (3) and the form of the coproduct in the
bialgebroid $(A\stac{B} A)^B$, the following identity holds in 
$(A\stac{B} A)^B\cten{R} (A\stac{B} A)^B$.
$$
\sum_{j\in J,m\in M_j,q} (x_q\stac{B} \gamma_j(x'_q))\cten{R}
({c^j}_m\stac{B} {c^{\prime j}}_m)=
\sum_{j\in J,m\in M_j,q}({c^j}_m\stac{B} {c^{\prime j}}_m)\cten{R}
(x_q\stac{B} \gamma_j(x'_q)).
$$
Hence, substituting the strong connection in Example \ref{ex.d2.strong} in
the definition of the ($k$-relative) Chern-Galois cycle of ${\bf e}$
in Theorem 
\ref{thm.galois.char}, its components come out as
$\tchg_0({\mathbf e})=\sum_{k\in K,l\in L,q} \varphi(v_k f_l {x}_q)
\omega( {x^{\prime}}_q e_l u_k)$ and, for $n\in \mathbb{N}$, 
\begin{eqnarray*}
\tchg_n({\mathbf e})
=\sum 
&\big( \bigotimes_{p=1}^{t-2} &\varphi( v_{k_p} f_{l_p} {c^{j_p}}_{m_p})
\omega(\gamma_{j_{p+1}}({c^{\prime j_p}}_{m_p}) e_{l_{p+1}} u_{k_{p+1}})\big)\\
&\otimes&\varphi(v_{k_{t-1}} f_{l_{t-1}} {c^{j_{t-1}}}_{m_{t-1}})
  \omega({c^{\prime j_{t-1}}}_{m_{t-1}} e_{l_{t}}
  u_{k_{t}}) \\
&\otimes&\varphi(v_{k_{t}} f_{l_{t}} {x}_q)
  \omega(\gamma_{j_{t}}({x^{\prime}}_q) e_{l_{t+1}}
  u_{k_{t+1}})\\ 
&\big( \bigotimes_{p=t+1}^{n+1} 
&\varphi( v_{k_p} f_{l_p}{c^{j_{p-1}}}_{m_{p-1}}) 
\omega(\gamma_{j_{p}}({c^{\prime j_{p-1}}}_{m_{p-1}}) e_{l_{p+1}} u_{k_{p+1}})\big)
\end{eqnarray*}
for any $t=1,\dots, n+1$, where 
the indices of $k$ and $l$ are understood modulo $n+1$
and the indices of $j$ and $m$ are  understood modulo $n$.
The summation is over the indices: $k_1,\dots , k_{n+1}\in K$, 
$l_1,\dots , l_{n+1}\in L$, $j_1,\dots ,j_{n}\in
J$, $m_1\in M_{j_1},\dots ,m_{n}\in M_{j_{n}}$ and $q$. 

Noting that $HC_0(B) = B / [B,B]$, and using the above explicit expression for
the components of the relative Chern-Galois cycle, one immediately finds that
the zeroth component of the relative Chern-Galois character has the following
simple form 
$\chg_0({\mathbf e})=\sum_{q} [\varphi({x^{\prime}}_q{x}_q)]
\in B / [B,B]$. 
\end{example}

The relative Chern-Galois character defined in Theorem~\ref{thm.galois.char}
is `relative', because it depends on $T$ and also on the choice of a strong
connection. The rest of this section is devoted to finding sufficient
conditions for $\chg_*^T$ to be independent of $\ell_T$. The basic idea
for this is to associate a finitely generated and projective $B$-module
to a given entwined extension and a comodule, and then to reveal
a connection between the relative
Chern-Galois
character and the relative Chern character of $B$.
 
Note that, given an entwined extension $B\subseteq A$ in
a right entwining structure $(A,\cC,\psi)_R$, and a left $\cC$-comodule 
 $W$, the cotensor product $A\Box_\cC W$ is a left $B$-module with the
 natural  
 product $b\stac{k} \sum_ia_i\stac{R} w_i\mapsto  \sum_iba_i\stac{R} w_i$. 
 We will study when this module is (relatively) projective. 
 
 \begin{theorem}\label{thm.proj} Let $(A,\cC,\psi)_R$ be a bijective right
 entwining structure over $R$ and let $B\subseteq A$ be an entwined extension.
 Let $T$ be a $k$-subalgebra of $B$. Suppose that
$A$ is a projective right $T$-module,
a flat left $T$-module and 
a faithfully flat right $B$-module.
If there exists a strong $T$-connection in $(A,\cC,\psi)_R$, then, for any
 left $\cC$-comodule $W$ that is finitely generated  
 as a left $R$-module,
 $
\Gamma =  A\Box_\cC W,
 $
 is a finitely generated and $T$-relative projective left $B$-module.
 \end{theorem}
\begin{proof}
 First we show that $\Gamma$ is a finitely generated left
 $B$-module. Since $A$ is a flat right $B$-module, there is 
an isomorphism of left $A$-modules,
 $
 A\stac{B}(A\Box_\cC W)\simeq (A\stac{B} A)\Box_\cC W.
 $
 In view of the 
 flatness of $A$ as a left or as a right $T$-module, the canonical map $\canA:
 A\stac{B}A\to A\stac{R}\cC$ is an isomorphism of right
 $\cC$-comodules by Corollary~\ref{cor.con}.
Since $-\Box_\cC W \, :\, \M^\cC\to\M_k$ is a covariant functor
(cf.\  \cite[21.3]{BrzWis:cor}), the map $\canA\Box_\cC W$ is 
 an isomorphism, so that the above  isomorphism can be extended to
 $$
 A\stac{B}(A\Box_\cC W)\simeq (A\stac{B} A)\Box_\cC W\simeq (A\stac{R}
 \cC)\Box_\cC W. 
 $$
 Finally, note that $(A\stac{R} \cC)\Box_\cC W\simeq A\stac{R} W$ 
 (cf.\ \cite[21.4-5]{BrzWis:cor}), hence there is an isomorphism of
 left $A$-modules 
 $$
 A\stac{B} \Gamma \simeq A\stac{R} W.
 $$ 
 Since $W$ is a finitely generated left $R$-module, $A\stac{R} W$ is a
 finitely generated 
 left $A$-module, hence also $A\stac{B} \Gamma$ is a finitely generated left 
 $A$-module. Since  
 $A$ is a faithfully
 flat right 
  $B$-module, \cite[Ch.\ 1\S 3 Prop.\ 11]{Bou:com} implies that $\Gamma$
  is a finitely generated left $B$-module.
  
Since $A$ is a projective right $T$-module by assumption and a
 $(\cC,T)$-projective right $B$-module by Theorem \ref{thm.Tstrconn}~(4), it
 follows that $A$ is a projective right $B$-module. As it is also a faithfully
 flat right $B$-module, the right regular $B$-module is
 a direct summand of  $A$ 
(cf.\ \cite[2.11.29]{Row:rin}). 
 In particular the right $T$-module $B$ is a direct summand of the flat
 $T$-module $A$, hence also $B$ is a flat right $T$-module. 

 Since all the assumptions of Theorem~\ref{thm.Tstrconn}~(3) are
 satisfied, there exists 
 a left $B$-module right $\cC$-comodule section $\sigma_T:A\to B\stac{T} A$ of
 the   product, and we can define the left $B$-module map
 $$
 \widetilde{\sigma} = \sigma_T\Box_\cC W : \Gamma = A\Box_\cC W\to
 (B\stac{T} A)\Box_\cC W \simeq B\stac{T} (A\Box_\cC W) = B\stac{T}\Gamma. 
 $$
 The last isomorphism follows by the flatness of $B$ as a right $T$-module.  
The map $\widetilde{\sigma}$ clearly is a section of the $B$-multiplication in
 $\Gamma$ (as $\sigma_T$ is a section of the multiplication map $B\stac{T}
 A\to A$), hence $\Gamma$ is a $T$-relative projective left $B$-module. 
  \end{proof}

 \begin{remark}\label{rem.dir.sum}
 In view of Corollary~\ref{cor.con2}, the assumption of Theorem~\ref{thm.proj} 
 that $A$ is a faithfully flat right $B$-module is
  satisfied provided $A$ is a faithfully flat right $T$-module.
 Furthermore, if the $\cC$-coaction on $A$ is given by a grouplike element 
 $e\in \cC$, then  $A$ is a faithfully flat right 
 $B$-module provided $B$ is a flat left $T$-module and a right
 $T$-module direct summand of $A$  
 by Corollary~\ref{cor:ff}. 
 (These assumptions hold e.g. in Example \ref{ex.d2.ch}).
Obviously, all these assumptions are 
 satisfied in case $T$ is equal to a ground field $k$.
 \end{remark}

 In relation to the Chern-Galois character, the case of a comodule $W$ that
  is not only finitely generated but also projective
 as a left $R$-module is of particular interest. In this case, 
 the $B$-module $\Gamma$ is projective under
 weakened assumptions on 
 the $T$-module $A$ and the explicit form of an idempotent 
 can be worked out.

Recall from \cite{Zim:pur} (cf.\ \cite[42.10]{BrzWis:cor}) that a
right $T$-module $M$ is  
 called a {\em locally projective module} if every finitely generated
 submodule of $M$ is projective, i.e.\ if, for any finitely generated 
 submodule $X$ of $M$, 
 there exist elements $x_1,\ldots, x_n\in M$ and $\xi_1, \ldots ,
 \xi_n \in M^*$  
 such that, for all $x\in X$, $x = \sum_{i=1}^n x_i\xi_i(x)$.  In
 particular, any 
  locally projective module is a flat module (cf.\ \cite[42.11]{BrzWis:cor}).
 The following theorem gives an explicit form of an idempotent for $\Gamma$ and
 thus asserts that $\Gamma$ is a finitely generated projective $B$-module.
\begin{theorem}\label{thm.idem}
Let $(A,\cC,\psi)_R$ be a bijective right entwining structure over $R$ and 
let $B\subseteq A$ be an entwined extension.
Let $T$ be a $k$-subalgebra of $B$ and suppose that the extension $B\subseteq
A$ splits as a $B$-$T$ bimodule. Suppose furthermore that $A$ is a locally
projective  
right $T$-module, a flat left $T$-module and that there exists a strong
$T$-connection  
 $\ell_T: \cC\to A\stac{T} A$ in  $(A,\cC,\psi)_R$.  Let $W$ be 
 a left $\cC$-comodule that is finitely generated and projective
 as a left $R$-module with a finite dual basis $\mathbf{w} = \{w_i\in
 W, \chi_i \in {}^* W\}_{i\in I}$, and let $e_{ij}\in \cC$ be as in
 Lemma~\ref{lemma.e}.  
 Let $X$ be a right $T$-submodule of $A$ finitely generated by
  $e_{ij}\suc 1_\nu$, $i,j\in I$, 
 where 
$ \ell_T(e_{ij}) = \sum_{\nu=1}^n e_{ij}\suc 1_\nu\stac{T} e_{ij}\suc 2_\nu$. 
Choose a finite set 
$\mathbf{x} = \{x_p\in A, \xi_p\in  \rhom T A T\}_{p\in P}$  
 such that, for all $x\in X$,  $x = \sum_{p\in P} x_p\xi_p(x)$, 
and define the family of right $\cC$-comodule maps
 $$
 \ell_p = (\xi_p\stac{T} A)\circ\ell_T : \cC\to A, \qquad p\in P.
 $$
 Choose a $B$-$T$ bimodule retraction $\phi: A\to B$ of the inclusion
 $B\subseteq A$ 
 and define the finite matrix of elements in $B$,  
 $$
 \mathbf{E} = (E_{(i,p),(j,q)})_{(i,p), (j,q)\in I\times P}, \qquad
 E_{(i,p),(j,q)} = \phi(\ell_p(e_{ij})x_q). 
 $$
 The matrix  $\mathbf{E}$ is an idempotent and the left $B$-module
 $\Gamma =  A\Box_\cC W$ is isomorphic to $ B^{(I\times P)}\mathbf{E}$.
 \end{theorem}
 
 \begin{remark}\label{rem.fflat}
 In view of Corollary~\ref{cor.con2}, in the case when $T$ is equal to the
 ground ring $k$, the assumptions of Theorem~\ref{thm.idem}  
 that $A$ is a locally projective $k$-module and $B$ is a direct summand
 in $A$ as a left $B$-module are satisfied provided $A$ is a
 faithfully flat and projective $k$-module. 
 Furthermore, by Proposition~\ref{prop:dirsum}, if the $\cC$-coaction in $A$
 is given by a grouplike element  
 $e\in \cC$ and $A$ is a locally projective $k$-module, then $B$ is a direct 
 summand in $A$ as a left $B$-module provided $B$ is a $k$-direct summand 
 of $A$. Obviously, all these assumptions are satisfied in case $k$ is a field.
 
 The assumptions of Theorem \ref{thm.idem} hold in Example \ref{ex.d2.ch}
 provided that $A$ is a locally projective $k$-module. 

 Note also that if, in addition to the assumptions of Theorem~\ref{thm.idem}, 
 the map \eqref{ups} is an
 epimorphism, then 
 $B\subseteq A$ is a $T$-flat extension.
 \end{remark}
The proof of the theorem uses the following lemma
formulated within
 the notation and assumptions of Theorem~\ref{thm.idem}.
\begin{lemma}\label{lemma.gamma}
For all $i\in I$ and $p\in P$, let
$
\gamma_{ip} = \sum_{j\in I} \ell_p(e_{ij})\stac{R} w_j \in A\stac{R} W.
$
Then
\begin{zlist}
\item $\gamma_{ip} \in \Gamma$;
\item $\sum_{j\in I, q\in P} E_{(i,p)(j,q)}\gamma_{jq} = \gamma_{ip}$.
\end{zlist}
\end{lemma}
\begin{proof}
(1) This is proven by the following explicit computation
\begin{eqnarray*}
(\roA \stac{R} W)(\gamma_{ip}) &=& \sum_{j\in I} \xi_p(e_{ij}\suc
1)e_{ij}\suc 2\sco 0 
\stac{R} e_{ij}\suc 2\sco 1\stac{R} w_j\\
&=& \sum_{j,k\in I} \xi_p(e_{ik}\suc 1)e_{ik}\suc 2
\stac{R} e_{kj}\stac{R} w_j\\
&=& \sum_{k\in I} \ell_p(e_{ik})\stac{R} \sum_{j\in I} e_{kj}\stac{R} w_j = 
(A\stac{R}\wro)(\gamma_{ip}),
\end{eqnarray*}
where the first equality follows by the left $T$-linearity of $\roA$ and
the second equality follows by the right $\cC$-colinearity of
the strong $T$-connection $\ell_T$ and by  
Lemma~\ref{lemma.e}~(3).

(2) By the  definition of the set $\mathbf{x}$, for all $c\in \cC$,
\begin{equation} \label{lemma.ell}
\sum_{p\in P}x_p\stac{T}\ell_p(c) = \ell_T(c).
\end{equation}
Compute
\begin{eqnarray*}
\sum_{j\in I, q\in P} E_{(i,p)(j,q)}\gamma_{jq} &=& \sum_{j,k\in I, q\in P} 
\phi(\ell_p(e_{ij})x_q)\ell_q(e_{jk})\stac{R}w_k\\
&=& \sum_{j,k\in I} 
\xi_p(e_{ij}\suc 1)\phi(e_{ij}\suc 2e_{jk}\suc 1)e_{jk}\suc 2\stac{R}w_k\\
&=& \sum_{j,k\in I} 
\xi_p(e_{ij}\suc 1)e_{ij}\suc 2e_{jk}\suc 1e_{jk}\suc 2\stac{R}w_k\\
&=&  
\sum_{k\in I} 
\xi_p(e_{ik}\suc 1)e_{ik}\suc 2\stac{R}w_k = \gamma_{ip},
\end{eqnarray*}
where the second equality follows by the $T$-$T$ bilinearity of $\phi$ and
\eqref{lemma.ell}. The third
equality follows by Lemma~\ref{lemma.e}~(3) combined with
Lemma~\ref{lemma.B} and the fact that $\phi$ restricted to $B$ is the
identity map. The penultimate equality is a direct consequence of the
definition of a strong $T$-connection (cf.\ proof of
Theorem~\ref{thm.Tstrconn}). 
\end{proof}

{\sl Proof of Theorem~\ref{thm.idem}.} We first show that $\mathbf{E}$
is an idempotent matrix. This is proven by the following direct
computation, for all $i,k\in I$ and $p,r\in P$, 
\begin{eqnarray*}
\sum_{j\in I, q\in P} E_{(i,p)(j,q)} E_{(j,q)(k,r)} &=&\sum_{j\in I,
q\in P} \phi(\ell_p(e_{ij})x_q)\phi (\ell_q(e_{jk})x_r) \\ 
&=& \sum_{j\in I} 
\xi_p(e_{ij}\suc 1)\phi(e_{ij}\suc 2e_{jk}\suc 1)\phi(e_{jk}\suc 2 x_r)\\
&=&  \sum_{j\in I} 
\xi_p(e_{ij}\suc 1)\phi(\phi(e_{ij}\suc 2e_{jk}\suc 1)e_{jk}\suc 2 x_r)\\
&=& \sum_{j\in I} 
\xi_p(e_{ij}\suc 1)\phi(e_{ij}\suc 2e_{jk}\suc 1e_{jk}\suc 2 x_r)
=  E_{(i,p)(k,r)},
\end{eqnarray*}
where the second equality follows by the $T$-$T$ bilinearity of $\phi$ and
\eqref{lemma.ell}, the third one is a  
consequence of the left $B$-linearity of $\phi$. The fourth equality
follows by 
Lemma~\ref{lemma.e}~(3) combined with Lemma~\ref{lemma.B} and the fact
that $\phi$ restricted to $B$ is the identity map, and the final
equality is 
a consequence of the definition of a strong $T$-connection. This concludes
the proof that 
$\mathbf{E}$ is an idempotent matrix. 

Consider the left $B$-module map
\begin{eqnarray*}
\Theta: B^{(I\times P)}\mathbf{E}&\to&\Gamma, \\
 \big(\sum_{i\in I,p\in P} b_{ip} E_{(i,p)(j,q)}\big)_{(j,q)\in
 I\times P} &\mapsto& \sum_{i,j\in I,p,q\in P} b_{ip}
 E_{(i,p)(j,q)}\gamma_{jq} = \sum_{i\in I,p\in P} b_{ip}\gamma_{ip}, 
\end{eqnarray*}
where the last equality follows by
Lemma~\ref{lemma.gamma}~(2). The map $\Theta$ has  its range
in  $\Gamma$ by
Lemma~\ref{lemma.gamma}~(1). We first show that $\Theta$ is an
injective map. Suppose that $\sum_{i\in I,p\in P}
b_{ip}\gamma_{ip}=0$. This implies that, for all $j\in I$, 
$$
0 = \sum_{i,k\in I, p\in P} b_{ip}\ell_p(e_{ik})\chi_j(w_k) =
\sum_{i,k\in I, p\in P} b_{ip}\ell_p(e_{ik}\chi_j(w_k)) = \sum_{i\in
I, p\in P} b_{ip}\ell_p(e_{ij}), 
$$
where we used that $\ell_p$ is a right $\cC$-comodule map, hence a
right $R$-module 
map, and Lemma~\ref{lemma.e}~(2). Thus, in particular, for all $j\in I$
and $q\in P$, $\sum_{i\in I, p\in P} b_{ip}\ell_p(e_{ij})x_q=0$,
hence, as $\phi$ is left $B$-linear, 
$$
0 = \sum_{i\in I, p\in P} b_{ip}\phi(\ell_p(e_{ij})x_q)=\sum_{i\in I,
p\in P} b_{ip}E_{(i,p)(j,q)}. 
$$
Therefore, $\Theta$ is a left $B$-module monomorphism. To prove that
$\Theta$ is  
an epimorphism we first take any $\sum_n a^n\stac{R} v^n\in \Gamma$
and compute 
\begin{eqnarray*}
\sum_{i,k\in I, p\in P, n} \!\!\!\!\!\!\!\!\!&
&\!\!\!\!a^n\chi_i(v^n) x_p\stac{T} 
\ell_p(e_{ik})\stac{R} w_k =
\sum_{i,k\in I, n} a^n\chi_i(v^n) \ell_T(e_{ik})\stac{R} w_k \\
&=& \!\!\!\!\! \sum_{i,k\in I, n}  \!\!\! a^n
\ell_T(\chi_i(v^n)e_{ik})\stac{R} w_k  
= 
\sum_{i\in I, n} a^n \ell_T((\chi_i(v^n)w_i)\sw{-1})\stac{R}
(\chi_i(v^n)w_i)\sw{0}\\ 
&=&  \sum_{n} a^n \ell_T(v^n\sw{-1})\stac{R} v^n\sw{0} 
=  \sum_{n} a^n\sco 0 \ell_T(a^n\sco{1})\stac{R} v^n,
\end{eqnarray*}
where the first equality follows by \eqref{lemma.ell}, the second
by the $R$-linearity of $\ell_T$ and the third by
Lemma~\ref{lemma.e}~(1). Then the dual basis property 
is used, and finally, the fact that $\sum_n a^n\stac{R} v^n$ is in the
cotensor product is employed. Recall from the proof of
Theorem~\ref{thm.Tstrconn} that the map $a\mapsto \sum a\sco 0
\ell_T(a\sco 1)\in B\stac{T} A$ is a section of the product $B\stac{T}
A\to A$, so that,  
applying $(\mu_A\stac{R} W)\circ(\phi\stac{T} A\stac{R} W)$ to the
equality just derived, we 
obtain
$$
\sum_n a^n\stac{R} v^n = \sum_{i\in I, p\in P, n} \phi(a^n\chi_i(v^n) x_p)
\gamma_{ip} .
$$
This shows that $\sum_n a^n\stac{R} v^n\in \im\Theta$, hence $\Theta$
is an epimorphism. Therefore,  $\Theta$ is a required  
isomorphism of left $B$-modules.
\endproof

With every finite
idempotent matrix $\mathbf{F} = (f_{ij})$ with entries from a $T$-ring
$B$ one associates a family of elements in the circular tensor product (cf.\ 
Section~\ref{sec.cyc})
\begin{equation}\label{chn}
\tch^T_{n}(\mathbf{F}) = \sum_{i_1,\ldots, i_{n+1}} f_{i_1i_2}\cten{T}
f_{i_2i_3}\cten{T}\cdots \stac{k} f_{i_{n}i_{n+1}}\cten{T}
f_{i_{n+1}i_1}\in B^{\cten{T}(n+1)}, 
\end{equation}
$n\in \mathbb{N}\cup\{0\}$. Since $\mathbf{F}$ is an idempotent and
$\tch^T_n(\mathbf{F})$ is invariant under the cyclic operator (modulo sign),
one easily finds that 
$\bigoplus_{l=0}^{2n} (-1)^{\lfloor l/2\rfloor}
\frac{l!}{\lfloor l/2\rfloor!}\tch^T_{l}(\mathbf{F})
$
 is a $2n$-cycle in the total relative
cycle homology complex (cf.\ the proof that $\chg^T_{2n}(\mathbf{e})$ is an
even cycle). The corresponding family of homology classes
\begin{equation}\label{chn.char}
\ch^T_{2n}(\mathbf{F}):=[\bigoplus_{l=0}^{2n} (-1)^{\lfloor l/2\rfloor}
\frac{l!}{\lfloor l/2\rfloor!}\tch^T_{l}(\mathbf{F})]\in HC_{2n}(B|T)
\end{equation}
 is called 
 a {\em $T$-relative Chern character of $B$}. Each of the elements 
$\tch_n^T({\bf F})$ defined by equation \eqref{chn} is called the
{\em $n$-th component of the  $T$-relative Chern cycle} associated to 
${\bf F}$. Note that the relative Chern character $\ch^T_{2n}(\mathbf{F})$ is
related to the Chern 
character $\ch_{2n}(\mathbf{F})$ (cf. \cite[Section~8.3]{Lod:cyc}), 
by $\ch^T_{2n}(\mathbf{F}) = \lambda_*\circ
\ch_{2n}(\mathbf{F})$, where $\lambda_*: HC_*(B)\to HC_*(B|T)$ is the canonical
surjection. Hence the relative Chern character defines a family of Abelian
 group morphisms $K_0(B)\to HC_*(B|T)$.

\begin{lemma}\label{lemma.ch} With the notation and assumptions of
Theorem~\ref{thm.idem}, the components of the $T$-relative 
Chern cycle associated to
$\mathbf{E}$ come out as $\tch_{0}^T(\mathbf{E}) =  
[\sum_{i\in I}\phi(e_{ii}\suc 2e_{ii}\suc 1)]$ and, for $n\in \mathbb{N}$,
\begin{eqnarray*}
\tch_{n}^T(\mathbf{E}) &=& \sum_{i_1,\ldots, i_{n+1}\in I} \!\!\!\!
e_{i_1i_2}\suc 2\ell_T(e_{i_2i_3})\ell_T(e_{i_3i_4})\cdots 
 \ell_T(e_{i_{n}i_{n+1}})\ell_T(e_{i_{n+1}i_1})e_{i_1i_2}\suc 1\\
& \!\!\!\! \!\!\!\!=& \!\! \!\!\!\!  \!\!\!\! \!\!\!\!\sum_{i_1,\ldots,
   i_{n+1}\in I}e_{i_1i_2}\suc 
2e_{i_2i_3}\suc 1 
\cten{T} e_{i_2i_3}\suc 2e_{i_3i_4}\suc 1\cten{T}\cdots 
\cten{T} e_{i_{n}i_{n+1}}\suc 2e_{i_{n+1}i_1}\suc 1
\cten{T} 
e_{i_{n+1}i_1}\suc 2e_{i_1i_2}\suc 1.
\end{eqnarray*}
\end{lemma}
\begin{proof} The explicit form of $\tch_n^T(\mathbf{E})$ is computed
directly from 
the definition of $\mathbf{E}$ and from equation \eqref{chn}. 
For the zeroth component,
$$
\tch_0^T(\mathbf{E}) = [\sum_{i\in I,
p\in P}  
\phi(\ell_{p}(e_{ii}) x_{p})] = [\sum_{i\in I,
p\in P}  
\phi(e_{ii}\suc 2 x_{p}\xi_p(e_{ii}\suc 1))] = [\sum_{i\in I}\phi(e_{ii}\suc
  2e_{ii}\suc 1)], 
$$
where the second equality follows by the definition of $A/[A,T]$ and the 
$T$-bilinearity of $\phi$.
With the help of 
\eqref{lemma.ell}, Lemma~\ref{lemma.e}~(3) and Lemma~\ref{lemma.B}, and
also using the fact that $\phi$ is a $T$-$T$ bimodule map and, if restricted to
$B$, it is the identity map, we can compute, for all $n>0$,
\begin{eqnarray*}
\tch_n^T(\mathbf{E})&=&\!\!\!\!\!\!\!\!\sum_{\stackrel{i_1,\ldots,
i_{n+1}\in I}{ p_1,\ldots, p_{n+1}\in P}}  
E_{(i_1,p_1)(i_2,p_2)}\cten{T}\cdots 
\cten{T} E_{(i_{n},p_{n})(i_{n+1},p_{n+1})} \cten{T} 
E_{(i_{n+1},p_{n+1})(i_1,p_1)}\\
&=&  \!\!\!\!\!\!\!\!\sum_{\stackrel{i_1,\ldots, i_{n+1}\in I}{
p_1,\ldots, p_{n+1}\in P}}  
\phi(\ell_{p_1}(e_{i_1i_2}) x_{p_2})\cten{T}\cdots  
\cten{T} \phi(\ell_{p_{n}}(e_{i_{n}i_{n+1}})x_{p_{n+1}}) \cten{T} 
\phi(\ell_{p_{n+1}}(e_{i_{n+1}i_{1}})x_{p_{1}})\\
&=&  \!\!\!\!\!\!\!\!\!\sum_{i_1,\ldots, i_{n+1}\in
I}\!\phi(e_{i_1i_2}\suc 2e_{i_2i_3}\suc 1) 
\cten{T} \cdots 
\cten{T} \phi(e_{i_{n}i_{n+1}}\suc 2e_{i_{n+1}i_1}\suc 1)\cten{T} 
\phi(e_{i_{n+1}i_1}\suc 2e_{i_1i_2}\suc 1)\\
&=&  \!\!\!\!\!\!\!\!\sum_{i_1,\ldots, i_{n+1}\in I} \!\!\!\!
e_{i_1i_2}\suc 2\ell_T(e_{i_2i_3})\ell_T(e_{i_3i_4})\cdots 
 \ell_T(e_{i_{n}i_{n+1}})\ell_T(e_{i_{n+1}i_1})e_{i_1i_2}\suc 1.
\end{eqnarray*}
This completes the proof.
\end{proof}

 If the map \eqref{ups} is an isomorphism, then 
 the components of the Chern character associated to $\mathbf{E}$ are
 equal to the components of the Chern-Galois character. This observation leads
 to the following
 \begin{theorem}
 Let $(A,\cC,\psi)_R$ be a bijective right entwining structure over $R$ and 
let $B\subseteq A$ be an entwined extension with a strong $T$-connection,
where $T$ is a $k$-subalgebra of $B$.  Suppose that the extension $B\subseteq
A$ splits as a $B$-$T$ bimodule, that $A$ is a locally projective 
right $T$-module, a flat left $T$-module and that the map  \eqref{ups}
is an epimorphism. Then the relative Chern-Galois character 
$\chg_{2n}^T: K_0(\cC)\to  HC_{2n}(B|T)$ does not depend on 
the choice of a  strong connection.
\end{theorem}
 \begin{proof}
 If $\Gamma = A\Box_\cC W$,  isomorphic $\cC$-comodules $W$ 
 lead to isomorphic $B$-modules $\Gamma$.  As cotensor product
 respects direct sums, the assignment $W\mapsto \Gamma$ descends
 to an Abelian map of Grothendieck groups. Thus, in view of
 Lemma~\ref{lemma.ch}, any component of the relative
 Chern-Galois character can be understood as a composition
 $$\xymatrix{
 \chg_{2n}^T: & K_0(\cC) \ar[r] &K_0(B)\ar[r]^{\ch_{2n}} 
 &HC_{2n}(B)\ar[r]^{\lambda_*} &HC_{2n}(B|T).}
 $$
The first map does not depend on the choice of $\ell_T$ as
 the definition of $\Gamma$ is independent of
 the choice of a strong connection. The second map is independent of
 $\ell_T$, since the Chern character
 is independent of the choice of an idempotent
 by  \cite[Theorem~8.3.4]{Lod:cyc}. The canonical epimorphism
 $\lambda_*: HC_{*}(B)\to HC_{*}(B|T)$ (cf. Section \ref{sec.cyc}) is obviously
 independent of $\ell_T$. 
  \end{proof}

\section*{Acknowledgements} 
The authors are grateful to Lars Kadison for his valuable comments.
The work on this paper started when the second author visited Budapest in
December  
2004. He would like to express his gratitude to the members of the Department
of 
Theoretical Physics at the Research Institute for Particle and Nuclear Physics
for very warm hospitality. 
The work of the first author is supported by the Hungarian Scientific 
Research Fund OTKA T 043 159, and of the second author by the EPSRC grant
GR/S01078/01.

\end{document}